\documentclass[lettersize,journal]{IEEEtran}
\usepackage{amsmath,amsfonts}
\usepackage{amssymb}
\usepackage{graphicx}
\usepackage[colorlinks=true, allcolors=blue]{hyperref}
\usepackage{bbm}
\usepackage{caption,subcaption}
\usepackage{verbatim}
\usepackage{cite,enumerate}
\usepackage{srcltx,cases}
\usepackage[mathscr]{eucal}
\usepackage{bbm}
\usepackage{multirow}

\newtheorem{theorem}{Theorem}

\newtheorem{corollary}{Corollary}
\newtheorem{remark}{\textit{Remark}}

\newtheorem{example}{Example}
\newtheorem{variant}{Variant}

\title{How Much Can Reconfigurable Intelligent Surfaces Augment Sky Visibility:\\
A Stochastic Geometry Approach}
\author{Junse Lee, Fran\c{c}ois Baccelli\thanks{Junse Lee is with School of AI Convergence, Sungshin Women's University, Seoul, South Korea (e-mail: junselee@sungshin.ac.kr).}
\thanks{Fran\c{c}ois Baccelli is with INRIA/ENS and Telecom Paris (e-mail francois.baccelli@ens.fr).}}

\begin{document}
\maketitle

\begin{abstract}
This paper uses the theory of point processes and stochastic geometry to quantify the sky visibility experienced by users located in an outdoor environment. The general idea is to represent the buildings of this environment as a stationary marked point process, where the points represent the building locations and the marks their heights. The point process framework is first used to characterize the distribution of the blockage angle, which limits the visibility of a typical user into the sky due to the obstruction by buildings. In the context of communications, this distribution is useful when users try to connect to the nodes of an aerial or non-terrestrial network in a Line-of-Sight way. Within this context, the point process framework can also be used to investigate the gain of connectivity obtained thanks to Reconfigurable Intelligent Surfaces.
Assuming that such surfaces are installed on the top of buildings to extend the user's sky visibility, this point process approach allows one to quantify the gain in visibility and hence the gain in connectivity obtained by the typical user. The distributional properties of visibility-related metrics are cross-validated by comparison to simulation results and 3GPP measurements.
\end{abstract}


\IEEEpeerreviewmaketitle
\section{Introduction}
\subsection{Motivation}

Through the 6G vision, the ITU has proposed recommendations on the target services that should be provided in the next-generation networks and on the core performance required in this context\cite{yazar20206g}. To accomplish this, academia and industry are actively considering introducing new network elements in next-generation communication systems to improve network performance. Recently, the 3rd Generation Partnership Project (3GPP) started exploring non-terrestrial networks (NTN) for expanding network services. NTN is a global term for any network involving flying objects such as high-altitude platform systems (HAPS), air-to-ground networks, and unmanned aerial vehicles (UAVs) \cite{lin20215g}. Future NTN networks will satisfy the stringent requirements of 6G while providing global coverage.

Understanding the connectivity of NTN entities to users on the ground is essential to optimize network planning and support seamless network services. In satellite networks, for example, the minimum elevation angle \cite{crowe1999model,li2002analytical} describes the visible region of the sky accessible by users on the ground, and the properties of evaluating this angle is essential to measure link performance and channel properties. 
Analyzing the sky visibility is essential to optimize the deployment of NTN nodes. However, previous frameworks proposed to evaluate this visibility have limitations and cannot be used to understand the interplay between network design parameters. 

Stochastic geometry provides a natural way to define and compute macroscopic properties of multi-user information theory communication channels. This tool is used to analyze the performance of wireless networks in terms of spatial averages. By using point processes to model the locations of network entities, stochastic geometry can provide valuable insights into the relation between network performance and network parameters. For example, stochastic geometry has been used to characterize the coverage and data rate of ad-hoc models \cite{baccelli2009stochastic,baccelli2006aloha,baccelli2009stochasticopp,NLee3}, cellular models \cite{andrews2011tractable,dhillon2012modeling,di2013average}, V2X \cite{tong2016stochastic,yi2019modeling}, and NTN networks \cite{chetlur2017downlink,banagar2020performance,okati2022nonhomogeneous,al2022next}. Stochastic geometry is also utilized to model blockage and shadowing effects by incorporating independent shadowing fields \cite{ilow1998analytic,blaszczyszyn2015wireless,bai2014analysis} and correlated shadowing fields \cite{baccelli2015correlated,zhang2015indoor,lee20163,lee2018effect}. In the present paper, we focus on modeling outdoor networks to measure users' Line-of-Sight (LOS) visibility into the sky by highlighting network environments (building locations and heights).

The paper is focused on the visibility and the coverage extension provided by reconfigurable intelligent surfaces (RISs). Several metrics are introduced to quantify this extension, like, e.g., the distribution of visibility angle gains due to RISs or the probability for the UE to be connected to a non-terrestrial node thanks to a RIS, given that it is not in the absence of RISs. Thanks to point process theory and stochastic geometry, these metrics are evaluated in closed forms. This in turn allows us to quantify the coverage extension and hence, in a sense, the economic value of RISs in this context.  Let us nevertheless stress a few obvious costs coming with these gains. First, RISs should be able to beamform to the satellite or HAPs and to the UE and should hence be adaptive and controlled and thus of a more expensive type. Second, detecting the RIS to be used as a relay (typically at the top of the building, limiting the sky visibility) will require either signaling and/or sensing. In order to estimate the cost-to-benefit ratio of such a RIS-enhanced architecture, the last two costs should also be evaluated. Even if the proposed model has the potential to quantify them, this will not be done here and will rather be left for future research.

\subsection{Related works}

There is a trend to use higher-frequency band signals as the sub-6 GHz spectrum band becomes saturated. Applications using high-frequency bands like mmWave communication can provide stable performance in the LOS environment, which is determined by network geometry. So, it is essential to understand the statistical behaviors of blockages and to identify the LOS paths.
In order to understand these statistical properties, many papers used stochastic geometry \cite{ilow1998analytic,blaszczyszyn2015wireless,bai2014analysis,baccelli2015correlated,zhang2015indoor,lee20163,lee2018effect}. Among them, \cite{bai2014analysis} used independent shadowing model to determine whether each path is LOS or NLOS by assuming that blockage and network nodes are independent. Novel frameworks to understand the correlated shadowing field were proposed for outdoor networks \cite{baccelli2015correlated}, indoor networks \cite{zhang2015indoor}, and in-building networks \cite{lee20163}. Also, \cite{lee2018effect} analyzed the effect on network performance when the shadowing effect is correlated. {\color{black}Recently, many papers have studied the LoS channel conditions in networks \cite{al2024line,bithas2024uav,roy2019new,tang2020tractable,koivumaki2021line,al2020probability,chen2023correlation}. The authors in \cite{al2024line,bithas2024uav,roy2019new,tang2020tractable,koivumaki2021line} proposed novel frameworks to investigate LoS channel conditions and performance analyses for outdoor networks. Especially, \cite{al2020probability,chen2023correlation} analyzed LoS probability between two points in an outdoor environment using the Boolean model of random shape theory.}

\subsection{Contributions}
This paper introduces a model to analyze the visibility into the sky experienced by a typical user in an outdoor network.

Let us stress that the basic model that we want to solve is fundamentally three-dimensional. There are two dimensions associated with the (assumed) planar structure of the city, as depicted in Figure \ref{fig:sys_model_2D}, where segments represent buildings seen from the sky. The third dimension, not depicted in this figure, is associated with the height of the buildings. This 3-D model determines the skyline that the users see. Assume that there are RISs installed at the top of all buildings
and that each user beamforms to the RIS at the top of the highest building in all directions of the outdoor plane around him/her. It should be clear that the model with random segments contains all needed information for evaluating the gain in sky visibility that each user gets by such a RIS extension.
Despite the fact that this is the ultimate aim of this line of research,
the analysis of the full 3-D problem will {\em not} be provided in the present paper.
The reason is that this question is too difficult at this stage.
{\color{black}We rather focus on investigating the visibility along the user's view direction by modeling and analyzing the simpler 2-D problem in closed form, where each user chooses a fixed or random direction in the plane and beamforms to the top of the building blocking his/her sight in this direction.}

In this 2-D model, which is depicted in Figure \ref{fig:sys_model_1D},
we represent the building locations as a homogeneous marked Poisson point process (PPP) along a one-dimensional line, and we assume that the building heights are i.i.d. random variables. Without loss of generality, we locate the typical user at the origin and define the visibility angle and the blockage angle into the sky observed by this the user. The main contributions on this 2-D model are summarized below:

\begin{itemize}
    \item {\color{black}We present the framework which is proposed to model the outdoor network using a marked point process. The location and mark of points represent the location and height of each building. We consider the situation where 1) the location of buildings follows a homogeneous PPP, and the height of buildings follows a Weibull distribution and its special case (Exponential distribution: M/M case). We provide two more cases for comparisons where 2) the location of buildings follows a homogeneous PPP and the height of the buildings is a fixed constant (M/D case), and 3) the distance between adjacent buildings is a fixed constant, and the height of buildings follows an exponential distribution (D/M case). }
    \item Based on these models, we analyze the user's visibility into the sky. We define the (positive) blockage angle as the maximum angle created by all buildings on the positive half-line and the visibility angle as $\frac{\pi}{2}-(\mbox{blockage angle})$. We derive the distributions of the user's blockage and visibility angles as functions of the network parameters. Then, we derive the joint distribution of the height and location of the building, which creates the blockage angle.
    \item We investigate the visibility enhancement obtained by leveraging RISs assumed to be installed at the top of each building. The user uses as relay the RIS on the building which creates its blockage angle. We consider two RIS operational modes, namely, the transmissive mode and the reflective mode. A transmissive RIS delivers the user's signal along the same direction while a reflective RIS delivers it along the opposite direction. We derive the distribution of the blockage angle at the RIS and investigate the visibility enhancement thanks to RISs. 
    \item We give integral expressions for metrics measuring how much RIS can enhance the visibility under the transmissive and reflective modes. We also compare the degree of visibility enhancement of the two RIS modes. 
\end{itemize}

The remainder of this paper is organized as follows. After presenting the network model and defining the blockage and the visibility angles in Section \ref{sec:Section 2}, we characterize the distribution of the two angles as functions of the network parameters and the joint distribution of height and location of the building which blocks the user's visibility into the sky in Section \ref{sec: Section3}. Then, we derive the visibility enhancement by RISs under the transmissive and reflective modes in Section \ref{sec: section4} and provide the numerical experiments to evaluate the visibility in Section \ref{sec: Section5} before concluding this paper in Section \ref{sec: Section6}.


\section{System Model}\label{sec:Section 2}

In this section, we present the network model for outdoor networks. In order to measure the user's visibility into the sky, we define two main geometric objects: the visibility and the blockage angle, which are the primary metrics experienced by the typical user.

\subsection{Network Model}

\begin{figure}[t]
     \centering
     \begin{subfigure}[b]{0.44\textwidth}
         \centering
         \includegraphics[scale=0.34]{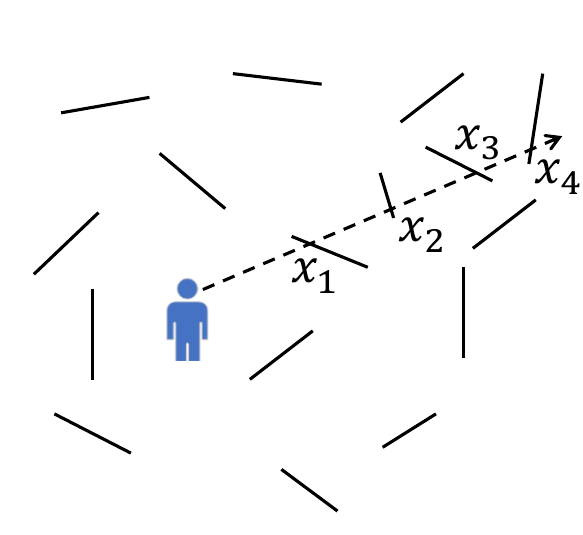}
         \caption{{\color{black}Top view of a three-dimensional outdoor network. The dotted arrow crossing the user represents the user's view direction. Line segments are terrestrial obstacles (buildings). }} 
         \label{fig:sys_model_2D} 
     \end{subfigure}
     \hfill
     \begin{subfigure}[b]{0.44\textwidth}
         \centering
         \includegraphics[scale=0.29]{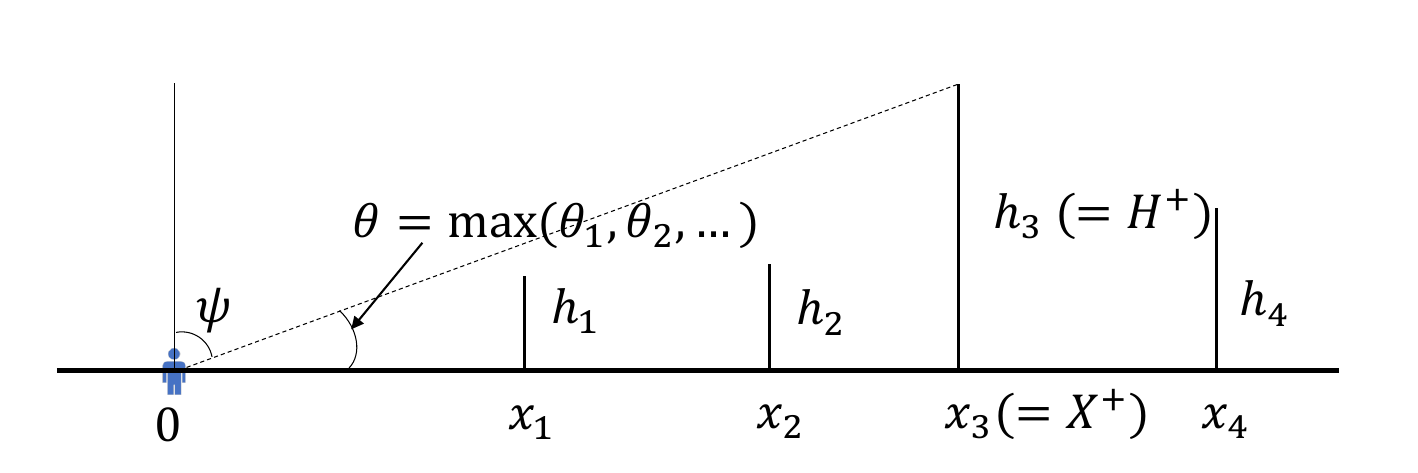}
         \caption{{\color{black}Two-dimensional network model along the user's view direction.}} 
         \label{fig:sys_model_1D} 
     \end{subfigure}
     \caption{Illustration of the outdoor network model.} 
     \label{fig:sys_model}
\end{figure}

 We consider a ground user who communicates with aerial network nodes. In order to understand the network performance experienced by such a user, a three-dimensional model must be used, with two dimensions being associated with the outdoor area of interest (represented as the Euclidean plane) and the third dimension of the building heights (represented as the positive half line). However, as already explained, such models are complex to analyze. So, in this paper, we focus on the user's view in a given direction to assess the user's visibility, which leads to a simpler (and tractable) two-dimensional model. The arrow passing through the user in Figure \ref{fig:sys_model_2D} represents the user's view direction, and this is modeled and analyzed as illustrated in Figure \ref{fig:sys_model_1D}. In Figure \ref{fig:sys_model_1D}, the locations and heights of the buildings on the line are modeled as a stationary marked point process. We denote the location of buildings as $\Phi = \{x_i\}$ which is a one-dimensional point process with intensity $\lambda$ and their associated marks as $H=\{h_i\}$ to represent their heights, with $\mu$ the inverse of the mean of $h_i$, respectively. The marks are assumed to be independent and identically distributed (i.i.d.) and independent of $\Phi$. We also define the parameter $\rho = \frac{\lambda}{\mu}$, which summarizes the nature of the outdoor network environment. Since $\rho$ is the product of the density of buildings and the mean height of buildings, a higher $\rho$ represents a denser outdoor environment or higher buildings.

\subsection{Visibility Angle and Blockage Angle}\label{subsec:sys_vis_angle_blk_angle}

Without loss of generality, we assume that the user is located at $(0,0)$. Let us denote by $\theta_i$ the angle with which the user sees the top of the building at $x_i$, as depicted in Figure \ref{fig:sys_model_1D}. In other words, the $\tan\theta_i$ is $\frac{h_i}{x_i}$. Since the user's visibility toward the sky along the positive direction is determined by the largest $\theta_i$, we define $\theta \triangleq \max \theta_i$ and call it the blockage angle observed by the user at $(0,0)$. Let $X^+$ and $H^+$ be the location and height of the building, which creates the blockage angle, $\theta$. So, $\tan \theta = \frac{H^+}{X^+}$. As $\theta$ decreases, the user visibility into the sky along the positive direction increases. We denote the visibility angle into the sky observed by the user along the positive direction by $\psi$, where $\psi = \frac{\pi}{2}-\theta$. 

In order to investigate the effect of the user's initial location and height, we further define $\theta_{x,h}$ and $\psi_{x,h}$ as the blockage angle and the visibility angle when the user's location is $(x,h)$ with $h\geq 0$. So, $\theta_{x,h}+\psi_{x,h} = \frac{\pi}{2}$ and $\theta = \theta_{0,0}$, $\psi = \psi_{0,0}$.

\subsection{Visibility Enhancement by RISs: Reflective and Transmissive Types}\label{subsec:vis_enh_RIS}
\label{ssec:ver}

In outdoor networks, network environments such as building locations and heights determine visibility into the sky, which in turn determines the connections between users and the aerial nodes. The simplest natural way to expand the user's visibility is to use a two-hop communication technique with the help of relay nodes. 

RISs can enhance the wireless communication environment and improve the signal quality \cite{elmossallamy2020reconfigurable,tang2020wireless}. The RIS technology uses new programmable subwavelength metamaterials to control electromagnetic waves intelligently. RISs are used for signal enhancement in scenarios such as blind spots of outdoor networks and indoor environments from outdoor base stations.

\begin{figure}[t]
\begin{center}
\includegraphics[scale=0.31]{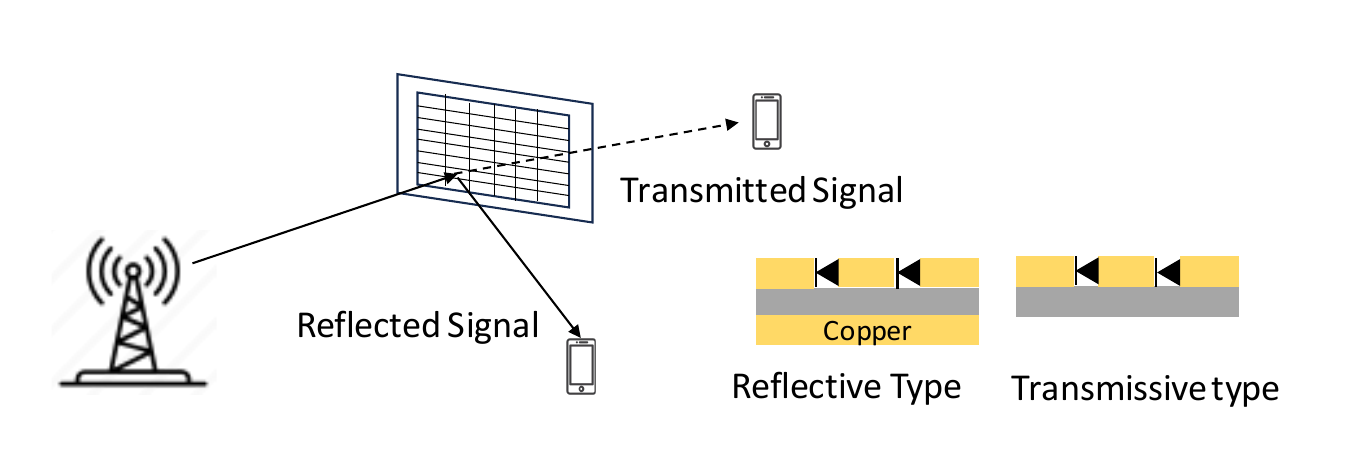}
\end{center}
\vspace{-0.5em} 
\caption{Operation modes and structures of RIS.} 
\label{fig:RIS} 
\end{figure}

RISs are categorized into three types: reflective, transmissive, and hybrid, as depicted in Figure \ref{fig:RIS}. Reflective RISs reflect signals toward the users on the same side of the base station (w.r.t. the RIS), while signals penetrate the transmissive type RIS so as to deliver signals to the users on the opposite side of the base station. Hybrid-type RISs enable a dual function of reflection and transmission. 

\begin{figure}[t]
     \centering
     \begin{subfigure}[b]{0.44\textwidth}
         \centering
         \includegraphics[scale=0.33]{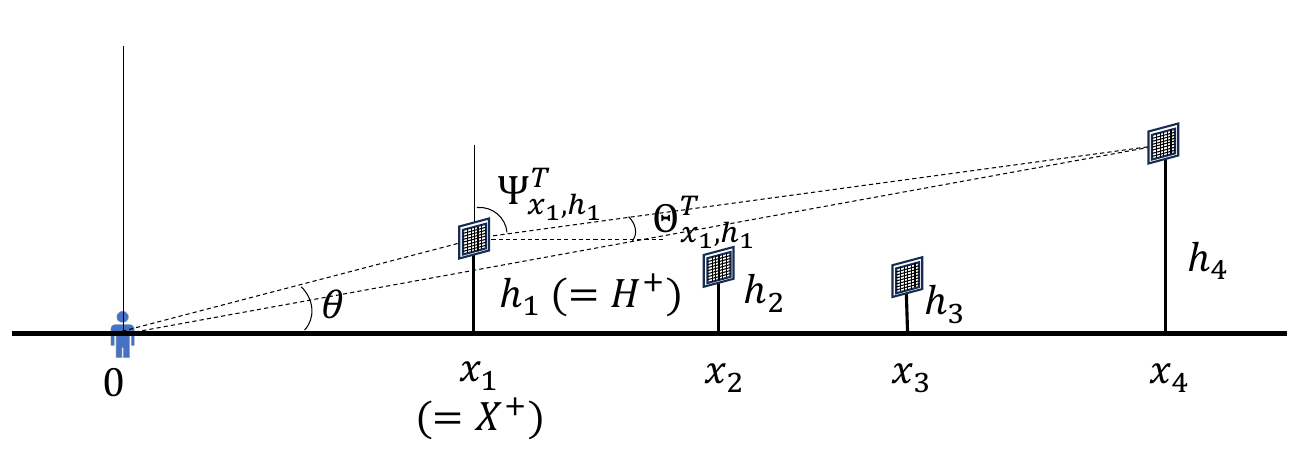}
         \caption{Transmissive mode} 
         \label{fig:sys_model1} 
     \end{subfigure}
     \vfill
     \begin{subfigure}[b]{0.44\textwidth}
         \centering
         \includegraphics[scale=0.33]{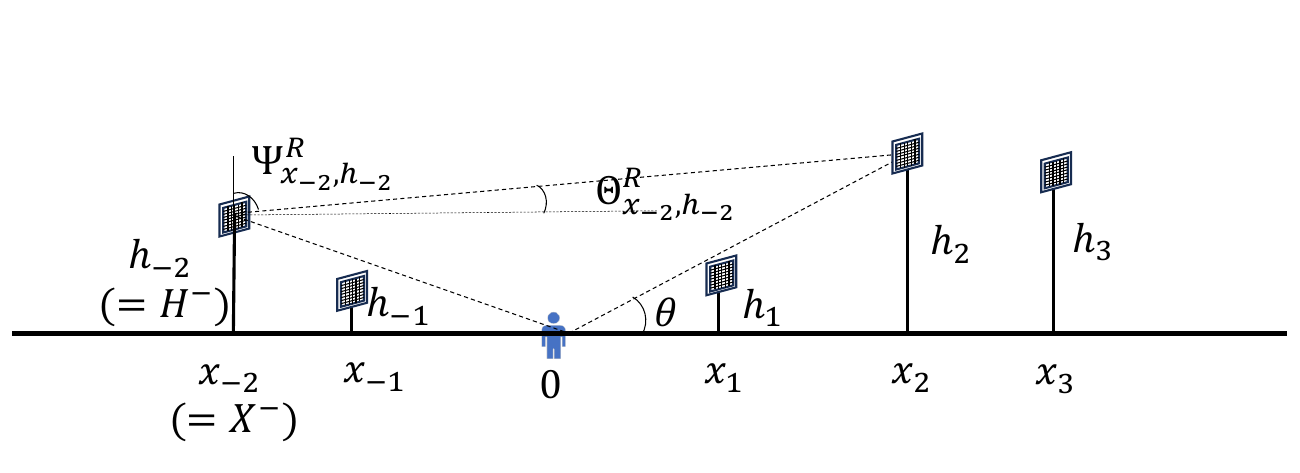}
         \caption{Reflective mode} 
         \label{fig:sys_model_left} 
     \end{subfigure}
     \caption{System model for visibility enhancement using RIS.}
     \label{fig:sys_model_all}
\end{figure}

In this paper, we consider RISs working as relay nodes in both reflective and transmissive modes to enhance visibility into the sky. We assume RISs are installed at the top of each building and that the user's signal is either transmitted through the RIS in the same direction or reflected by the RIS in the opposite direction. With the aid of RISs, the ground user can reach non-terrestrial nodes that are not visible at their position on the ground in a two-hop manner. {\color{black}RISs are capable of active beam steering toward a desired direction. Since we focus on visibility along the user's view direction, we assume that the beam can be steered in either direction of transmission or reflection, as illustrated in Figure \ref{fig:sys_model_all}.}

Figure \ref{fig:sys_model1} illustrates a scenario where a transmissive RIS enhances the visibility of the user at $(0,0)$. The user's signal arrives at a transmissive RIS installed on the roof of the building, which creates $\theta$. This signal penetrates the RIS and is relayed to aerial nodes. In order to analyze the visibility enhancement using transmissive RISs, we define the blockage angle and the visibility angle when the user's signal is retransmitted by a transmissive RIS installed at $(X^+,H^+)$ given that $(X^+,H^+)=(x,h)$ as ${{\Theta}}_{x,h}^T$ and ${{\Psi}}_{x,h}^T$, respectively\footnote{We use capital letters for these RIS centric angles to distinguish them from the user-centric angles, for which we use lower cases.}.

The visibility enhancement by a reflective RIS is illustrated in Figure \ref{fig:sys_model_left},  when the user is located at $(0,0)$. We define $X^-$ and $H^-$ as the location and height of the building which blocks the visibility toward the sky along the negative direction. To take advantage of the reflective property, the user transmits the signal along the negative direction, and the signal is then reflected along the positive direction at the reflective RIS installed on the top of the building at $X^- < 0$ with a height $H^->0$. We denote by ${\Theta}^R_{x,h}$ and ${\Psi}^R_{x,h}$ the (RIS-centric)  blockage angle and visibility angles in the positive direction, namely the blockage and visibility angles viewed from $(X^-,H^-)$ given that $(X^-,H^-)=(x,h)$, respectively. 

We further define the following symbols to quantify how much RISs improve the sky visibility. Let $\Theta^T$ and $\Theta^R$ be the blockage angles under the transmissive mode and the reflective mode obtained by deconditioning $\Theta_{x,h}^T$ and $\Theta_{x,h}^R$ with respect to $x$ and $h$, respectively. In other words, ${\Theta}^T = \mathbb{E}_{x,h}[{\Theta}^T_{x,h}]$ and ${\Theta}^R=\mathbb{E}_{x,h}[\Theta^R_{x,h}]$. Similarly, we also define $\Psi^T$ and $\Psi^R$ as the unconditioned visibility angles under the transmissive and reflective modes with respect to $x$ and $h$, respectively. So, $\Theta^T+\Psi^T=\Theta^R+\Psi^R=\frac{\pi}{2}$.

Table I summarizes the geometric parameters used in our network model.

{ \tiny
\begin{table}
\begin{center}
\begin{tabular}[]{|c|l|} 
 \hline
 \textbf{Symbol} & \textbf{Definition}  \\ 
  \hline
 $\lambda$  & Density of buildings  \\ 
 \hline
  $\mu$ & Inverse of the mean height of buildings\\
 \hline
   \multirow{2}{*}{$\rho (= \frac{\lambda}{\mu})$}  & Product of the density of buildings \\
   & and the mean height of buildings   \\ 
 \hline
 \multirow{2}{*}{$\theta$, $\psi$}  & Blockage angle and visibility angle \\
 & observed by the user at $(0,0)$  \\ 
  \hline
 \multirow{2}{*}{$\theta_{x,h}$, $\psi_{x,h}$} & Blockage angle and visibility angle \\
 & observed by the user at $(x,h)$  \\
  \hline
 \multirow{2}{*}{$X^+$, $H^+$} & Location and height of the blocking building\\
 &along the user's view direction  \\ 
 \hline
 \multirow{2}{*}{$X^-$, $H^-$}  & Location and height of the blocking building \\
 &along the opposite of the user's view direction  \\ 
 \hline
  \multirow{3}{*}{$\Theta^T_{x,h}$, $\Psi^T_{x,h}$}  & Blockage angle and visibility angle observed \\
  &by the transmissive RIS at $(x,h)$\\
  &given that $X^+=x$, $H^+=h$ \\ 
  \hline
 \multirow{3}{*}{$\Theta^R_{x,h}$, $\Psi^R_{x,h}$} & Blockage angle and visibility angle observed \\
 &by the reflective RIS at $(x,h)$\\
 &given that $X^-=x$, $H^-=h$   \\ 
   \hline
 \multirow{2}{*}{$\Theta^T$, $\Psi^T$} & Blockage angle and visibility angle \\
 &observed by the transmissive RIS  \\ 
 \hline
  \multirow{2}{*}{$\Theta^R$, $\Psi^R$} & Blockage angle and visibility angle \\
  &observed by the reflective RIS  \\ 
  \hline
\end{tabular}
\end{center}
\vspace{-1.5em}
\caption{List of symbols for our network model.}
\end{table}
}

\section{Visibility Analysis}\label{sec: Section3}

\subsection{Max Shot Noise Analysis of the Visibility Angle}
We first analyze the distribution of $\theta$ and $\psi$ depicted in Figure \ref{fig:sys_model}, which characterize the user's visibility into the sky along the positive direction. Let $f(x_i,h_i)=\frac{h_i}{x_i}=\tan \theta_i $, and $f(X^+,H^+)=\frac{H^+}{X^+}=\tan \theta = \tan \underset{i}{\max} \theta_i = \underset{i}{\max} \tan \theta_i = \underset{i}{\max} f(x_i,h_i)$. 

{\color{black} We derive the distributions of $\theta$ when the location of buildings is a realization of a homogeneous PPP with intensity $\lambda$ under the assumption that 1) the height of buildings follows a Weibull distribution with shape parameter $k$ and scale parameter $\frac{1}{\mu}$. We also consider a special case by plugging $k=1$, so that the Weibull distribution reduces to an exponential distribution with parameter $\mu$ (M/M case), 2) the height of buildings is fixed to a given constant $\frac{1}{\mu}$ (M/D case). Further, we derive the distribution of $\theta$ under the assumption that 3) the distance between adjacent buildings is fixed as $\frac{1}{\lambda}$ and the height of buildings follows an exponential distribution with parameter $\mu$ (D/M case). In these cases, $\Phi$ is stationary so that the distribution of the visibility angles seen by the user at $0$ is the same as that seen by a user at $s\in\mathbb{R}$ for all $s$.  }

{\color{black}
In the next theorem, we present the distribution of $\tan\theta$ when the height of buildings follows a Weibull distribution with scale parameter $k$ and shape parameter $\frac{1}{\mu}$.

\theorem The cumulative distribution function (CDF) of $\tan\theta$ is 
\begin{equation}
    \mathbb{P}[\tan\theta \leq t] = \exp\left({-\frac{\rho}{t}\Gamma\left(1+\frac{1}{k}\right)}\right), ~~\mbox{$t\geq 0$.}
\end{equation}
where $\Gamma(\cdot)$ is the Gamma function.
\begin{IEEEproof}
See Appendix \ref{appen1}.
\end{IEEEproof}

The following corollary gives the distribution of $\theta$.
\corollary The CDF of $\theta$ is 
\begin{equation}
    \mathbb{P}[\theta \leq \phi] = \exp\left({-\frac{\rho}{ \tan \phi}\Gamma\left(1+\frac{1}{k}\right)}\right)
\end{equation}
and the probability density function (PDF) of $\theta$ is
\begin{equation}
f_{\theta}(\phi) = \exp\left({-\frac{\rho}{ \tan \phi}\Gamma\left(1+\frac{1}{k}\right)}\right) \frac{\rho\Gamma\left(1+\frac{1}{k}\right)}{ \sin^2\phi},
\end{equation}
for $\phi\in[0,\frac{\pi}{2})$, respectively.

\begin{IEEEproof}
{\color{black}The CDF and PDF of $\theta$ are obtained by applying a change of variables in (1) and differentiating (2) with respect to $\phi$, respectively. }
\end{IEEEproof}

Now, we analyze the M/M case by plugging $k=1$ in Theroem 1 and Corollary 1.
}

\example \label{theo1} In the M/M case, the CDF of the $\tan\theta$ is 
\begin{equation}
    \mathbb{P}[\tan\theta \leq t] = e^{-\frac{\rho}{t}}, ~~\mbox{$t\geq 0$,}
\end{equation}
which is a Fr\'echet distribution with shape parameter $\alpha=1$ and scale parameter $s = \rho$.

From the characteristics of the Fr\'echet distribution, we get that the mean (and hence the second moment) of the $\tan{\theta}$ is infinite.

\example\label{cor1} In the M/M case, the CDF of $\theta$ is 
\begin{equation}\label{eq:corollary1-1}
    \mathbb{P}[\theta \leq \phi] = \exp\left({-\frac{\rho}{ \tan \phi}}\right),
\end{equation}
and the probability density function (PDF) of $\theta$ is
\begin{equation}\label{eq:corollary1-2}
f_{\theta}(\phi) = \exp\left({-\frac{\rho}{ \tan \phi}}\right) \frac{\rho}{ \sin^2\phi},
\end{equation}
for $\phi\in[0,\frac{\pi}{2})$, respectively.

Example \ref{cor1} shows that the CDF of ${\theta}$ is a decreasing function of $\rho$. From the relation of stochastic dominance, we can conclude that increasing the density and height of buildings is likely to increase the blockage angle, $\theta$, which is in line with intuition.

\begin{figure}[t]
\begin{center}
\includegraphics[scale=0.20]{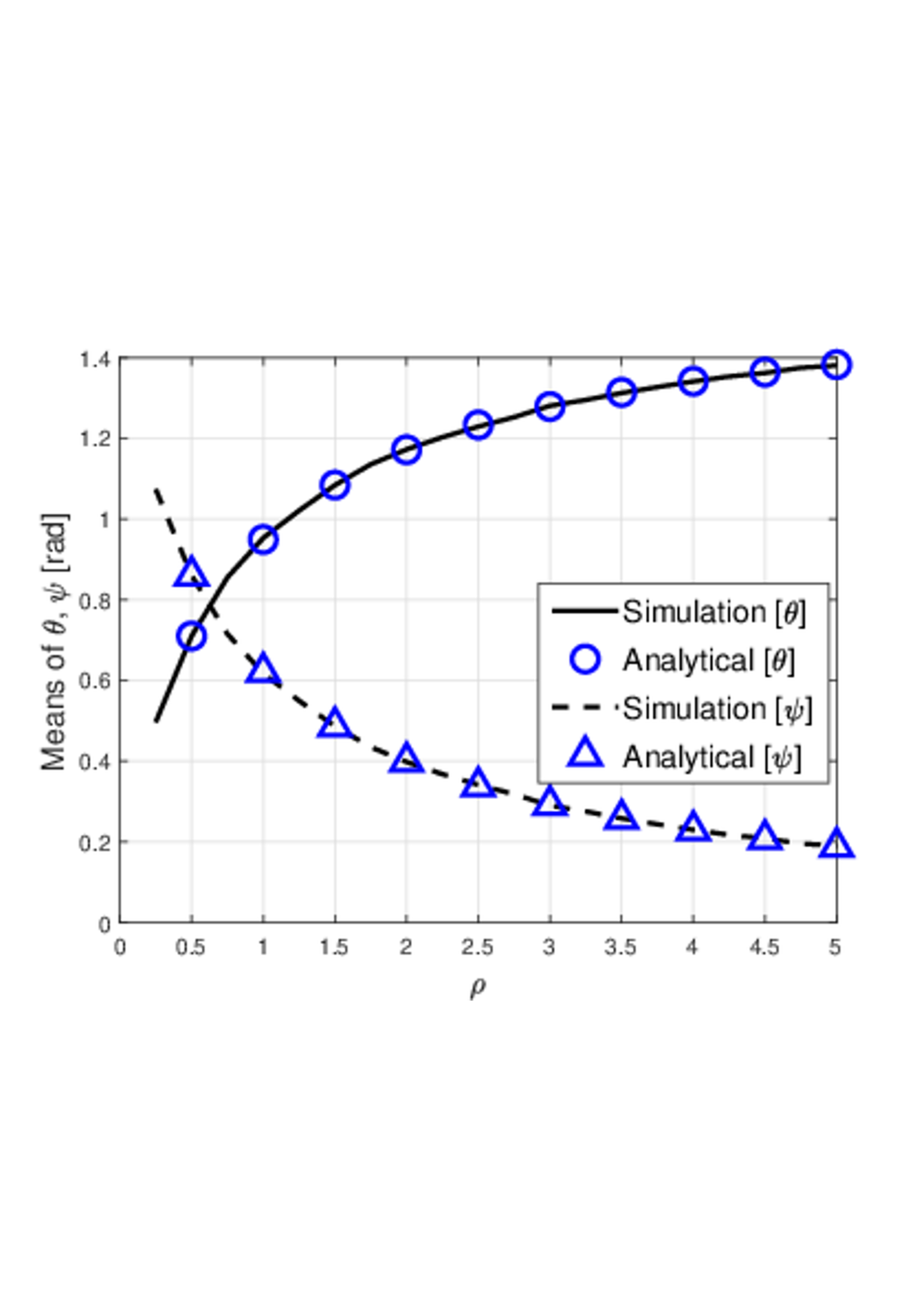}
\end{center}
\vspace{-0.5em} 
\caption{Means of ${\theta}$ and ${\psi}$ in the M/M case with respect to $\rho$. } 
\label{fig:mean_angles} 
\end{figure}

The PDF of $\theta$ provides the mean blockage angle as
\begin{align}\label{eq:mean_blockage_angle}
&\mathbb{E}[{\theta}] =    \int_{0}^{\frac{\pi}{2}} \phi f_{{\theta}}(\phi) d\phi 
=\Bigg[\exp\left({-\frac{\rho\Gamma\left(1+\frac{1}{k}\right)}{ \tan \phi}}\right)\phi\nonumber\\
+& \frac{i}{2}\Bigg( e^{i \rho\Gamma\left(1+\frac{1}{k}\right)} Ei\left(-\frac{\rho\Gamma\left(1+\frac{1}{k}\right)}{2}\left(i+ \cot \frac{\phi}{2}\right)^2 \tan \frac{\phi}{2}\right) \nonumber\\
+&  e^{-i \rho\Gamma\left(1+\frac{1}{k}\right)} Ei\left(\frac{\rho\Gamma\left(1+\frac{1}{k}\right)}{2}\cot \frac{\phi}{2} \left(i+\tan \frac{\phi}{2}\right)^2\right)   \Bigg)\Bigg]_0^{\frac{\pi}{2}} \nonumber\\
=& \frac{1}{2\sqrt{\pi}}G_{2,4}^{3,2}\left(
\begin{array}{c}
\frac{1}{2},1\\
\frac{1}{2},\frac{1}{2},1,0
\end{array}\middle\vert
\frac{\left(\rho\Gamma\left(1+\frac{1}{k}\right)\right)^2}{4}
\right),
\end{align}
where $Ei(z) = -\int_{-z}^{\infty} \frac{e^{-t}}{t} \mathrm{d}t$ is the exponential integral function \cite{spanier1987exponential} and $G^{\cdot,\cdot}_{\cdot,\cdot}(\cdot|\cdot)$ is the Meijer G-function \cite{beals2013meijer}. The mean of ${\psi}$ is obtained as $\mathbb{E}[{\psi}]=\pi/2-\mathbb{E}[{\theta}]$. Figure \ref{fig:mean_angles} illustrates the mean of ${\theta}$ and ${\psi}$ in function of $\rho$ in the M/M case $(k=1)$. {\color{black}As $\rho$, the product of the density of buildings and mean height of buildings, increases, the mean blockage angle also increases. This is aligned with the fact that the user's visibility into the sky is more limited when the density of the buildings is denser and the heights of buildings are larger.} Unlike the mean of the $\tan{\theta}$, which is equal to infinity, the mean of ${\theta}$ is bounded by $\frac{\pi}{2}$.

\corollary The CDF and the PDF of $\psi$ are
\begin{align}
\mathbb{P}[\psi \leq \phi ] &= 1-\exp\left(-\rho\tan \phi\right),\nonumber\\
f_{\psi}(\phi)&=\frac{\rho}{\cos^2\phi} \exp\left(-\rho\tan \phi\right),\nonumber
\end{align}
in the M/M case.
\begin{IEEEproof}
The distribution of $\psi$ is obtained from the relation $ \psi = \frac{\pi}{2}-\phi$.
\end{IEEEproof}

We now consider the M/D and D/M cases to be compared with the M/M case in Theorem \ref{theo1}. The first variant considers the situation where all buildings' heights are constant and equal (M/D), and the second one is when all buildings are located periodically (D/M). For a fair comparison, we assume that the buildings' height is $\frac{1}{\mu}$ in the M/D case and that the distance between adjacent buildings is $\frac{1}{\lambda}$ in the D/M case.

\variant \label{ex1} In the M/D case, the CDF of $\tan{\theta}$ is  
\begin{equation}\label{eq:MDCDF}
    \mathbb{P}[\tan{\theta} \leq t] \nonumber \\ = \exp\left(-\frac{\rho}{ t}\right),
\end{equation}
and the CDF and PDF of $\theta$ are 
\begin{align}
    \mathbb{P}[{\theta} \leq \phi] &= \exp\left(-\frac{\rho}{ \tan \phi}\right),\\
    f_{{\theta}}(\phi)&=\exp\left(-\frac{\rho}{ \tan \phi}\right)\frac{\rho}{ \sin^2 \phi},
\end{align}
respectively, for $\phi\in[0,\frac{\pi}{2})$.
\begin{IEEEproof}
{\color{black}  When the building's height is constant, ${\theta}$ is achieved by the first building located at $x_1$. 
Since $x_1$ follows an exponential distribution with parameter $\lambda$, the CDF of $\theta$ is obtained as \eqref{eq:MDCDF}. The PDF of $\theta$ is derived by differentiating it.}
\end{IEEEproof}

From Example \ref{theo1} and Variant \ref{ex1}, we can observe that the distributions of $\theta$ are the same whether the height of buildings follows an exponential distribution (M/M) or is constant (M/D).

\variant\label{ex2} In the D/M case, we suppose that $x_1$ is uniformly located in $[0,\frac{1}{\lambda})$, and that the distance between adjacent buildings is $\frac{1}{\lambda}$. This makes the point process $\sum\delta_{x_i}$ stationary \cite{baccelli2020random}. If the building's heights follow an exponential distribution with parameter $\mu$, the CDF of $\tan{\theta}$ becomes
\begin{align}\label{eq6}
    &\mathbb{P}[\tan {\theta} \leq t]=\nonumber\\
    &\lambda\int_{0}^{\frac{1}{\lambda}}\prod_{i=1}^{\infty}\left(1-\exp\left(-\mu\left(u+\frac{i-1}{\lambda}\right)t\right)\right)\mathrm{d}u,
\end{align}
and the CDF of ${\theta}$ becomes
\begin{align}\label{eq7}
    &\mathbb{P}[{\theta} \leq \phi]=\nonumber\\
    &\lambda\int_{0}^{\frac{1}{\lambda}}\prod_{i=1}^{\infty}\left(1-\exp\left(-\mu\left(u+\frac{i-1}{\lambda}\right)\tan \phi\right)\right)\mathrm{d}u,
\end{align}
for $\phi\in[0,\frac{\pi}{2})$.

\begin{IEEEproof}
See Appendix \ref{appenvar2}.
\end{IEEEproof}

\begin{figure}[t]
\begin{center}
\includegraphics[scale=0.20]{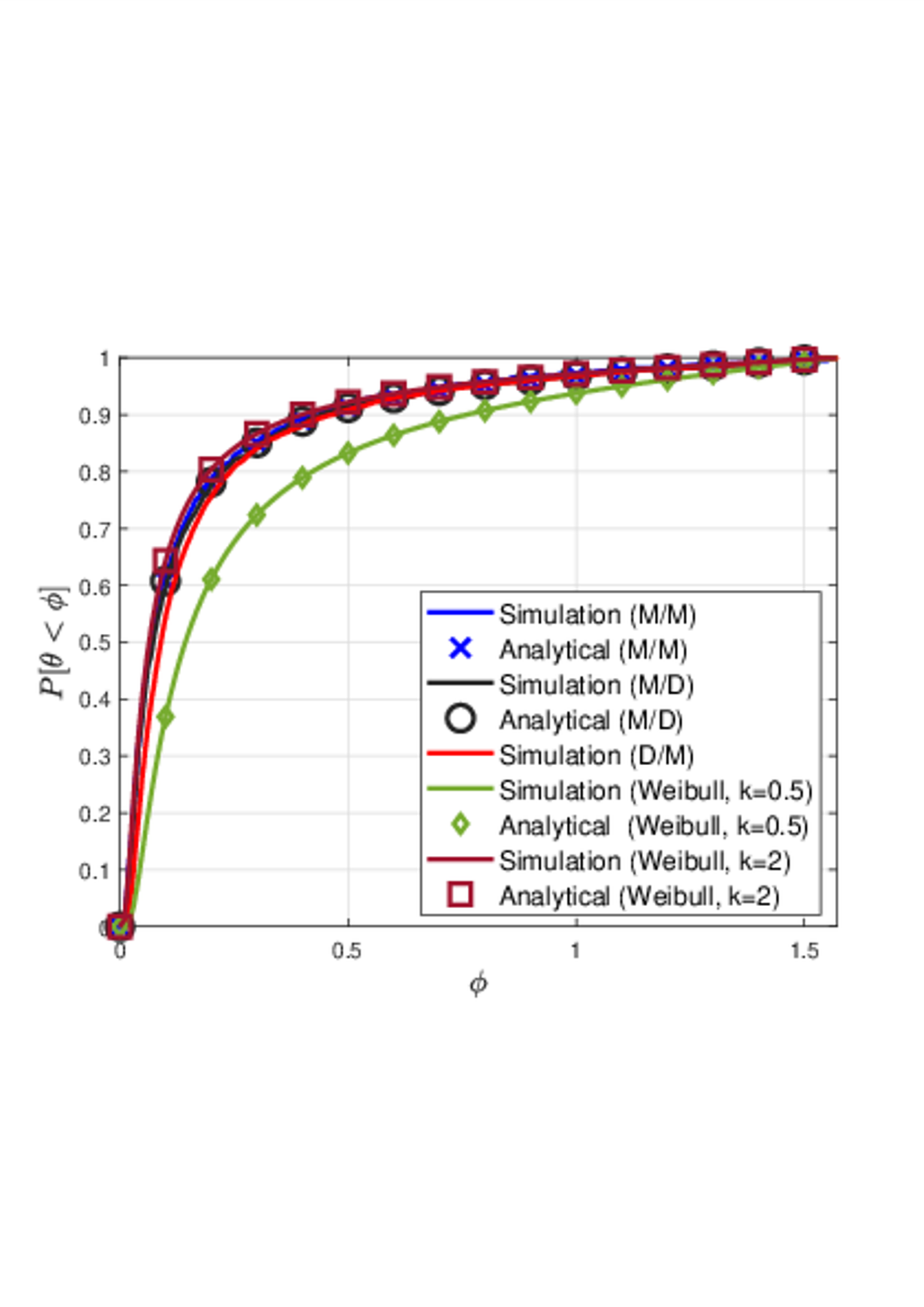}
\end{center}
\vspace{-0.5em} 
\caption{CDF of $\theta$ under $\lambda = 0.1$ and $\mu = 2$.} 
\label{fig:blockage_angle_dist}
\end{figure}

The distribution of $\theta$ and that of its tangent are illustrated in Figure \ref{fig:blockage_angle_dist}. We set $\lambda = 0.1$ and $\mu = 2$, which gives $\rho = 0.05$. 

So far, we have analyzed the distribution of $\theta$ when the user's position is $(0,0)$. The following theorem provides the distribution for the blockage angle under the M/M case for general user's location.  

\theorem \label{theo3} In the M/M case, the CDF of $\tan{\theta}_{x,h}$ is 
\begin{align}\label{eq:tanxh}
    \mathbb{P}[\tan {{\theta}}_{x,h}\leq t]&=\exp\left(-\frac{\rho }{ t} e^{-\mu h}\right),
\end{align}
for $t\geq 0$, whereas the CDF and the PDF of ${\theta}_{x,h}$ are
\begin{align}\label{eq:CDF_theta+}
    \mathbb{P}[{\theta}_{x,h}\leq \phi] &= \exp\left(-\frac{\rho e^{-\mu h}}{ \tan \phi} \right),\\
    f_{{\theta}_{{x},h}}(\phi) & =  \exp\left(-\frac{\rho e^{-\mu h}}{ \tan \phi} \right)\frac{\rho e^{- \mu h}}{\sin^2{\phi}},
\end{align}
for $\phi \in [0,\frac{\pi}{2})$, respectively.
\begin{IEEEproof}
See Appendix \ref{appentheo3}.
\end{IEEEproof}

Theorem \ref{theo3} shows that the distribution of ${\theta}_{x,h}$ is independent of $x$, which comes from the stationarity of the homogeneous PPP. So, the visibility enhancement induced by the user's location only comes from his/her height of the user. The visibility enhancement can be checked by comparing the stochastic ordering of the CDF of $\theta$ and that of $\theta_{x,h}$, i.e., $\mathbb{P}[{\theta}_{x,h}\leq \phi] \geq \mathbb{P}[{\theta}_{0,0}\leq \phi]=\mathbb{P}[\theta \leq \phi] $ for all $\phi$. Further, Theorem \ref{theo3} implies that the gain of visibility obtained by increasing the height ${h}$ is the same as the gain obtained by decreasing the density of the building by $\lambda e^{-\mu h}$ when comparing with Corollary \ref{cor1}. This effect comes from the thinning theorem of the PPP \cite{baccelli2009stochastic} and the memoryless property of the exponential random variable. Since the probability that the height of each building is larger than $h$ is $e^{-\mu h}$, the visibility of the user located at $h$ is affected by buildings which form a homogeneous PPP with intensity $\lambda e^{-\mu h}$. This explains why Theorem \ref{theo3} boils down to Theorem \ref{theo1} by plugging in $h=0$.

\begin{figure}[t]
\begin{center}
\includegraphics[scale=0.20]{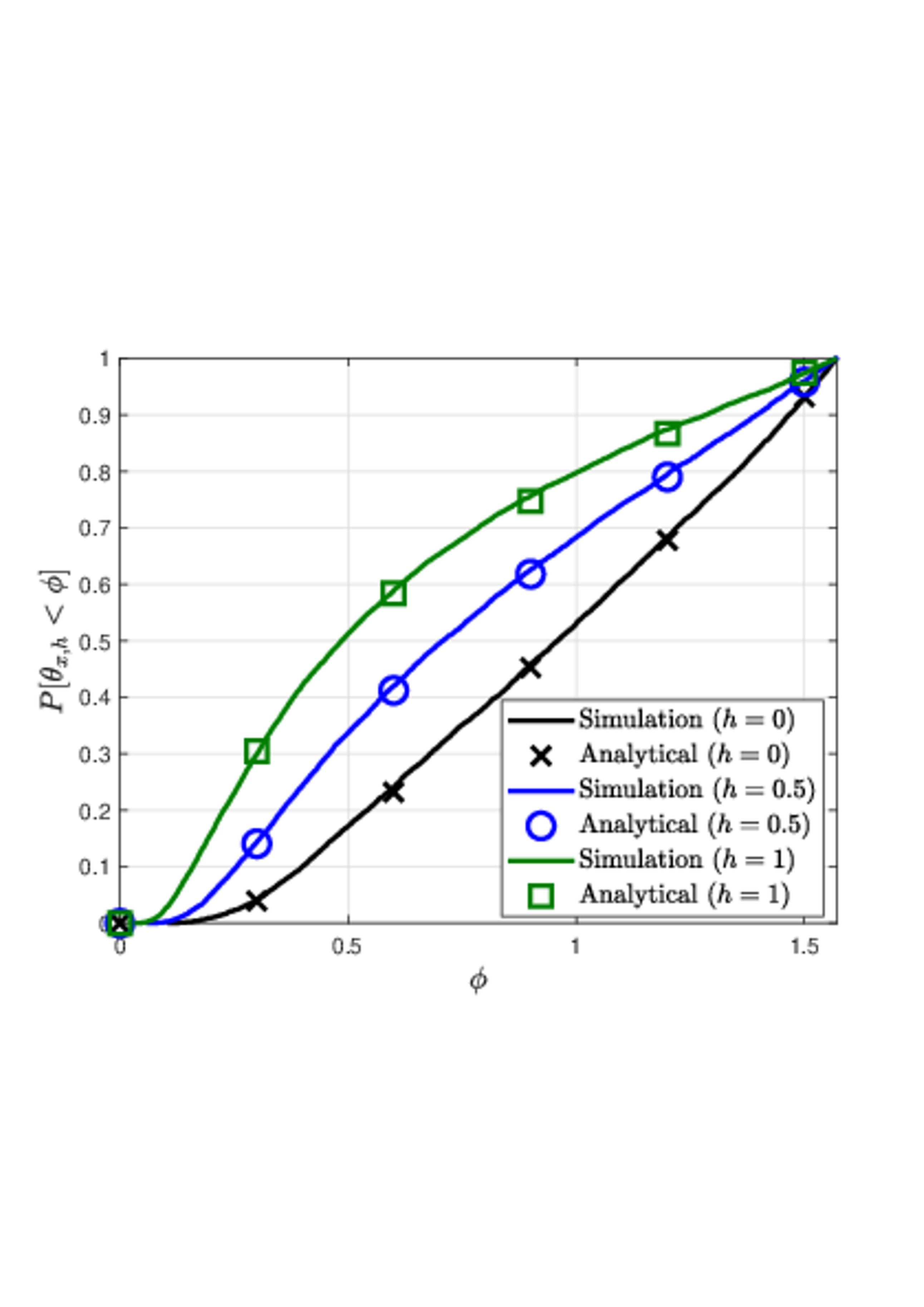}
\end{center}
\vspace{-0.5em} 
\caption{CDFs of ${\theta}_{x,h}$ when $\lambda = 1$ and $\mu = 1$.} 
\label{fig:vis_gain_right} 
\end{figure}

Figure \ref{fig:vis_gain_right} illustrates the stochastic ordering comparison of ${\theta}_{x,h}$ for $\lambda = 1$ and $\mu = 1$. We can observe that, for all $\phi$, $\mathbb{P}[{\theta}_{x,h_1}\leq \phi] \geq \mathbb{P}[ {\theta}_{x,h_2}\leq \phi] $ when $h_1\geq h_2$, which shows the visibility enhancement obtained from the increase of the user's initial height.

\begin{figure}[t]
\begin{center}
\includegraphics[scale=0.20]{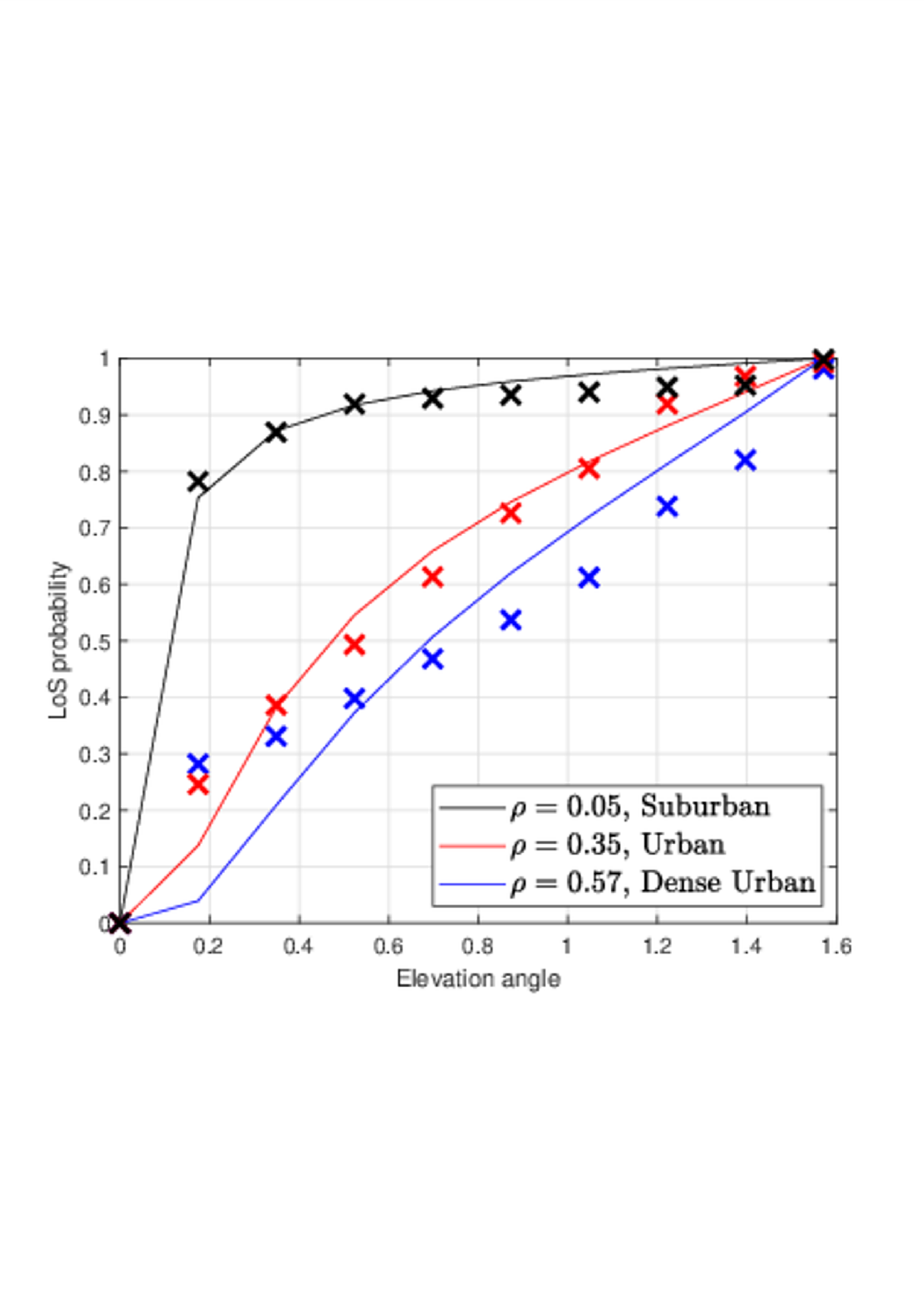}
\end{center}
\vspace{-0.5em} 
\caption{Comparison of LOS probabilities based on our model and 3GPP data \cite{3gpp38811} } 
\label{fig:response1_3_2} 
\end{figure}

\example(Connection with the minimum elevation angle) Theorem \ref{theo3} provides an analytical expression for the minimum elevation angle, which is an important metric to measure the link quality of satellite networks. For a user located at $(x,h)$ who observes the sky with an angle of $\zeta$ with respect to the ground, the probability that his/her sight is blocked by obstacles is $\mathbb{P}[\theta_{x,h}\leq \zeta]=\exp\left(-\rho e^{-\mu h}\cot{\zeta}\right)$. To the best of our knowledge, the analytical expressions for the minimum elevation angle capturing how network environment parameters affect that angle are new. {\color{black}In Figure \ref{fig:response1_3_2}, we compare our model (M/M, M/D)-based LOS probability and the empirical measurements of 3GPP \cite{3gpp38811}. The markers indicate 3GPP data, and we used parameters $\rho = 0.05, 0.35, 0.57$ and $h=0$ to fit the suburban, urban, and dense urban scenarios. This figure shows that results from our M/M and M/D models are well-matched with 3GPP data.
}

\corollary The CDF and the PDF of ${\psi}_{x,h}$ are
\begin{align}\label{eq:CDF_psi+}
   \mathbb{P}[{\psi}_{x,h}\leq \phi] &= 1-\exp\left(-{\rho e^{-\mu h} \tan \phi }\right),\\
    f_{{\psi}_{x,h}}(\phi) & =  \exp\left(-{\rho e^{-\mu h} \tan \phi }\right)\frac{\rho e^{- \mu h}}{\cos^2{\phi}}.
\end{align}
The CDF and the PDF of ${\psi}_{x,h}$ are derived from Theorem \ref{theo3} and the relation ${\psi}_{x,h}=\frac{\pi}{2}-{\theta}_{x,h}$.

The mean of $ {\theta}_{x,h}$ is
\begin{align}\label{eq:meanthetaxh}
&\mathbb{E}[{\theta}_{x,h}]  = \int_{0}^{\frac{\pi}{2}} \phi f_{{\theta}_{x,h}}(\phi) \mathrm{d}\phi  =\Bigg[e^{-\frac{e^{-h\mu}\rho}{\tan \phi}}\phi \nonumber \\
&+ \frac{i e^{-i e^{-h\mu}\rho}}{2}\Bigg( e^{2 i e^{-h\mu}\rho} Ei\Bigg(-\frac{\rho e^{-h \mu \rho}}{2} \nonumber\\
&\times\left(i+ \cot \frac{\phi}{2}\right)^2\tan \frac{\phi}{2}\Bigg)\nonumber\\
&- Ei\Bigg(\frac{e^{-h \mu}\rho}{2}\cot \frac{\phi}{2} \left(i+\tan \frac{\phi}{2}\right)^2\Bigg)   \Bigg)\Bigg]_0^{\frac{\pi}{2}}  \nonumber\\
&= \frac{1}{2\sqrt{\pi}}G_{2,4}^{3,2}\left(
\begin{array}{c}
\frac{1}{2},1\\
\frac{1}{2},\frac{1}{2},1,0
\end{array}\middle\vert
\frac{(\rho e^{-h\mu})^2 }{4}
\right),
\end{align}
and the mean of ${\psi}_{x,h}$ is
\begin{align}\label{eq:barpsi}
&\mathbb{E}[{\psi}_{x,h}]  = \frac{\pi}{2}-\mathbb{E}[{\theta}_{x,h}].
\end{align}
\begin{figure}[t]
\begin{center}
\includegraphics[scale=0.20]{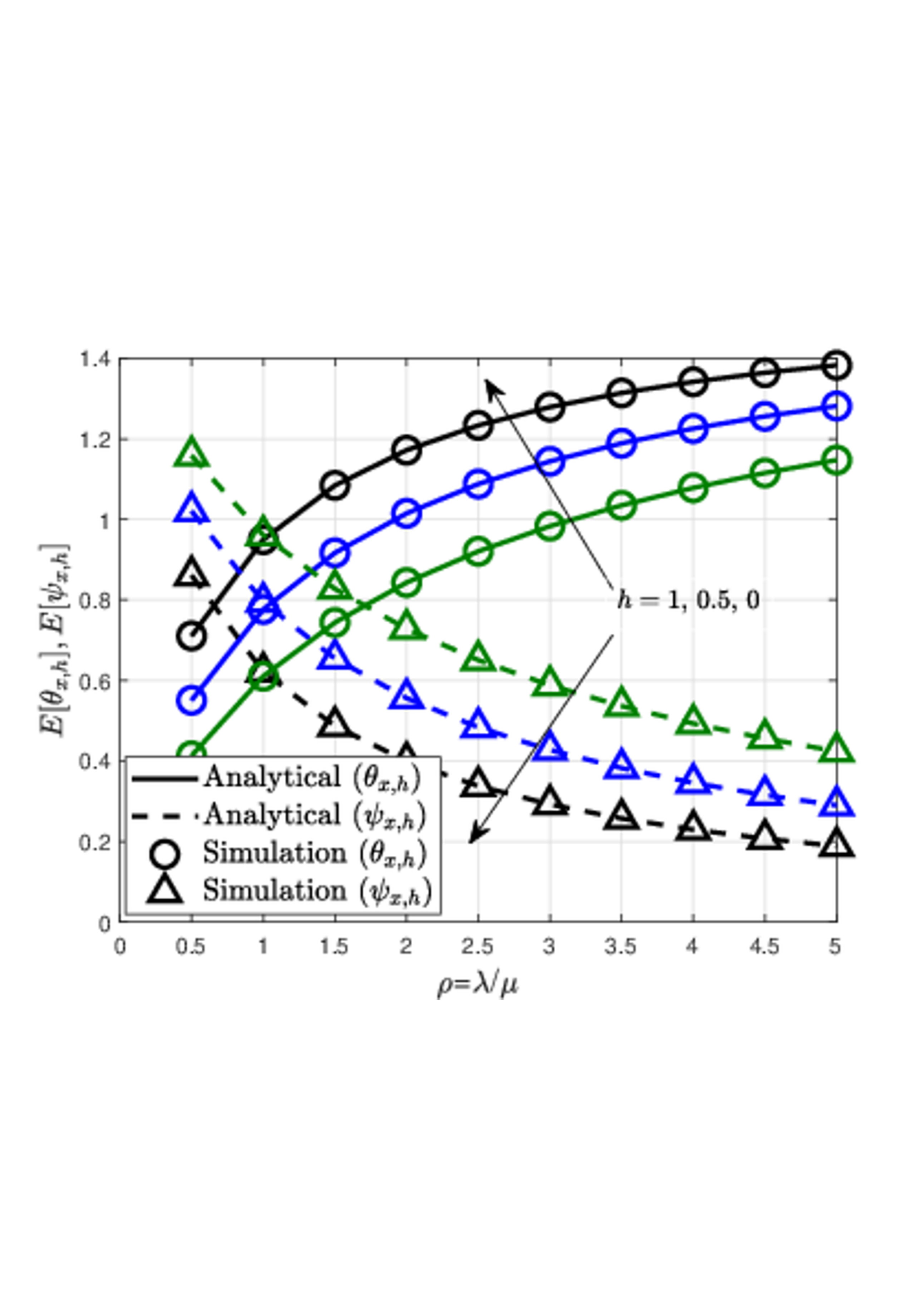}
\end{center}
\vspace{-0.5em} 
\caption{Means of $\theta_{x,h}$ and $\psi_{x,h}$ with different $h$.} 
\label{fig:mean_enh} 
\end{figure}

Figure \ref{fig:mean_enh} shows that \eqref{eq:meanthetaxh} and \eqref{eq:barpsi} match well with the means of $\theta_{x,h}$ and $\psi_{x,h}$ obtained from simulations with different $h$. This figure corroborates our intuition that the user's initial height affects its visibility. As $h$ decreases, the blockage angle at $(x,h)$, $\theta_{x,h}$, increases and the visibility angle at $(x,h)$, $\psi_{x,h}$, decreases.

\subsection{Joint Distribution}

In this subsection, we focus on the M/M case and the case where the user is at $(0,0)$.

The following theorem provides the joint distribution of $X^+$ and $H^+$, which are the location and the height of the building which achieves $\theta$ as defined in Section \ref{subsec:sys_vis_angle_blk_angle}.

\theorem \label{theo2} The joint density of $(X^+,H^+)$ at $(x,h)$ is
\begin{equation}\label{eq:joint}
j(x,h)=\lambda \mu e^{-\mu h - \frac{\lambda x}{\mu h}}.
\end{equation}

\begin{IEEEproof}
See Appendix \ref{appentheo2}.
\end{IEEEproof}

In the next corollary, we provide the marginal distributions of $H^+$ and $X^+$. 
\corollary\label{cor2} The density of $H^+$ at $h$ is 
\begin{equation}\label{eq:joint_mar_x}
g(h)=\mu^2 h e^{-\mu h},
\end{equation}
which is a Gamma distribution of parameters $(2,\frac 1 \mu)$
and the density of $X^+$ at $x$ is
\begin{equation}\label{eq:joint_mar_h}
k(x)=2\lambda\sqrt{\lambda x} K_1(2\sqrt{\lambda x}),
\end{equation}
where $K_n(\cdot)$ is the modified Bessel function of the second kind \cite{abramowitz1968handbook}. 
\begin{IEEEproof}
The marginal distributions are obtained by integrating \eqref{eq:joint} w.r.t $x$ and $h$, respectively.
\end{IEEEproof}

We can verify that \eqref{eq:joint} is a PDF since integrating $g(h)$ w.r.t. $h$ is 1. Further, we see from \eqref{eq:joint} and Corollary \ref{cor2} that $H^+$ and $X^+$ are not independent.

\begin{figure}[t]
\begin{center}
\includegraphics[scale=0.20]{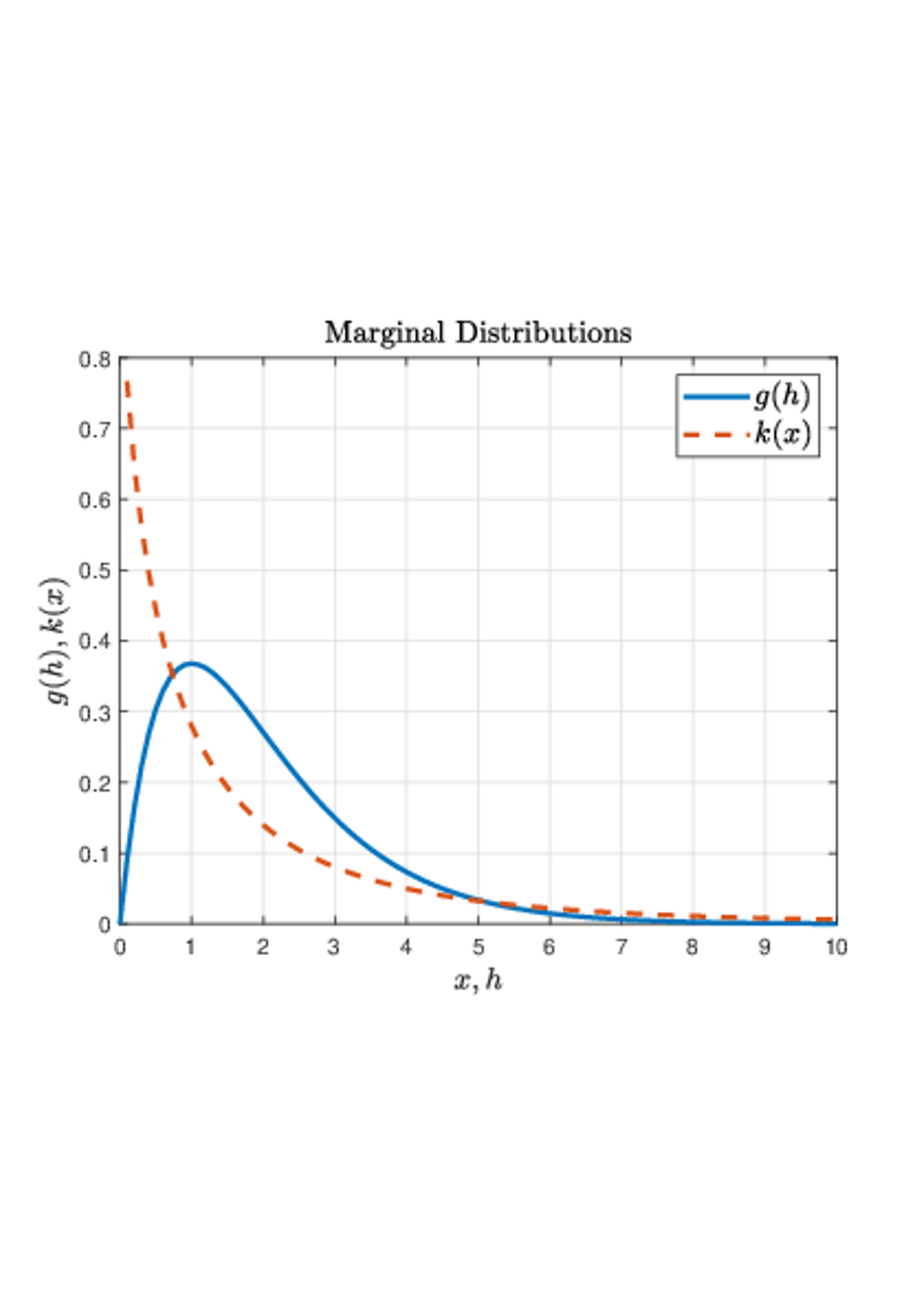}
\end{center}
\vspace{-0.5em} 
\caption{Marginal distributions of $X^+$ and $H^+$ under $\lambda, \mu = 1$.} 
\label{fig:joint} 
\end{figure}

Figure \ref{fig:joint} depicts the marginal distributions of $H^+$ and $X^+$ when $\lambda, \mu = 1$. By differentiating \eqref{eq:joint_mar_x}, we can find that $g(h)$ has a maximum value at $\frac{1}{\mu}$ and this is also illustrated in Figure \ref{fig:joint}.

\corollary\label{cor3} The means of $H^+$ and $X^+$ are $\frac{2}{\mu}$ and $\frac{2}{\lambda}$, respectively.
\begin{IEEEproof}
    \begin{align}
        \mathbb{E}[H^+] &= \int_{0}^{\infty} h g(h) \mathrm{d}h= \frac{2}{\mu}.\nonumber\\
        \mathbb{E}[X^+]&=\int_{0}^{\infty} x k(x) \mathrm{d}x = \frac{2}{\lambda}.\nonumber
    \end{align}\end{IEEEproof}
The average height of the building constituting the blockage angle is $\frac{1}{\mu}$ in the M/D case, while it is $\frac{2}{\mu}$ in the M/M case even though the distributions of $\theta$ are identical for the M/M and M/D cases as shown in the previous subsection.

\begin{remark}
By symmetry, the joint distribution of $(-X^-,H^-)$ coincides with that of
$(X^+,H^+)$.
\end{remark}

\section{Visibility Enhancement by RISs}\label{sec: section4}

A natural way to increase the user's two-hop visibility into the sky is to install additional relays to communicate in a two-hop manner. In this subsection, we consider a network of RISs installed on the roofs of all buildings. These are used as relays for the users' signals toward aerial network nodes to improve the users' visibility. The analysis is conducted under the M/M case. We further assume that the user is located at $(0,0)$. We investigate the two-hop visibility enhancement induced by RISs under the transmissive and the reflective mode.

\subsection{Transmissive Mode}\label{subsec:tr}

Figure \ref{fig:sys_model1} illustrates a scenario where a transmissive RIS enhances the user's visibility in a two-hop manner. The user's signal arrives at a transmissive RIS installed on the roof of the building, which creates $\theta$. This signal penetrates the RIS and is delivered to aerial nodes on the other side.

The next theorem provides the distribution of $\tan{\Theta}_{x,h}^T$ (defined in Subsection \ref{ssec:ver}).

\theorem\label{theo:vis_en_tr_tan} The conditional CDF of $\tan{\Theta}^T$ given that $(X^+,H^+)=(x,h)$ is
\begin{align}\label{eq:vis_en_tr_tan}
     &\mathbb{P}[\tan {\Theta}_{x,h}^T \leq t] = \mathbb{P}[\tan {\Theta}^T \leq t |X^+=x, H^+=h]  \nonumber\\
    &=
    \begin{cases}
    \exp\left(-\rho\left[\frac{1}{t}-\frac{x}{h}\right]e^{-\mu h}\right), & \text{for $0<t\leq\frac{h}{x}$,}\\
    1, & \text{for $t>\frac{h}{x}$.}
    \end{cases}
\end{align}

\begin{IEEEproof}
See Appendix \ref{appentheo_vis_en_tr_tan}.
\end{IEEEproof}

When $x=0$, \eqref{eq:vis_en_tr_tan} reduces to \eqref{eq:tanxh} in Theorem \ref{theo3}, since the visibility angle observed by transmissive RIS installed at $(0,h)$ and that observed by the user at $(0,h)$ are identical.

Since we assume that $(X^+,H^+) = (x,h)$, we have $\frac{h_i}{x_i}\leq \frac{h}{x}$ for all $i$. This condition introduces the discontinuity of the PDF of $\tan\Theta_{x,h}^T$ identified in the next corollary.

\corollary The conditional density of $\tan {\Theta}^T$ given $(X^+,h^+)=(x,h)$ is equal to
\begin{align}\label{eq:trans_tan}
    &f_{\tan {\Theta}_{x,h}^T}(t) =\nonumber\\
    &
    \begin{cases}
    \frac{\rho}{t^2}  \exp\left(-\mu h -e^{-\mu h}\rho\left(\frac{1}{t}-\frac{x}{h}\right)\right), & \text{for $0<t\leq\frac{h}{x}$} \\
    0. & \text{otherwise.}
    \end{cases}
\end{align}
\begin{IEEEproof}
$f_{\tan {\Theta}^T_{x,h}}(t)$ is obtained by differentiating \eqref{eq:vis_en_tr_tan} with respect to $t$.
\end{IEEEproof}
\corollary The $k$-th order moment of $\tan {\Theta}^T$ given $(X^+,H^+)=(x,h)$ is 
\begin{align}
    &\mathbb{E}[(\tan {\Theta}_{x,h}^T)^k]=c(x,h) \alpha^{k-1} \Gamma\left(-k+1, \frac{\alpha x}{h}\right),
\end{align}
where $c(x,h) = \frac{\lambda}{\mu}e^{\frac{\lambda x}{\mu h}e^{-\mu h}-\mu h}$, $\alpha = \frac{\lambda}{\mu}e^{-\mu h}$, and $\Gamma(s,x)=\int_{x}^{\infty}t^{s-1}e^{-t}\mathrm{d}t$ is the upper incomplete Gamma function.
\begin{IEEEproof}
    The $k$-th order moment is
    \begin{align}
       &\mathbb{E}[(\tan {\Theta}^T_{x,h})^k]=\mathbb{E}[(\tan {\Theta}^T)^k|X^+=x, H^+=h]\nonumber\\
       &= c(x,h)\int_{0}^{\frac{h}{x}} t^{k-2} \exp\left(-\frac{\lambda}{\mu t} e^{-\mu h}\right) \mathrm{d}t\nonumber\\
        & = c(x,h) \int_{\frac{x}{h}}^{\infty} \frac{1}{v^{k-2}} \exp\left({-\frac{\lambda}{\mu } e^{-\mu h}v}\right) \frac{1}{v^2} \mathrm{d}v \nonumber\\
        & = c(x,h) \int_{\frac{x}{h}}^{\infty} v^{-k} e^{-\alpha v} \mathrm{d}v \nonumber\\
        & = c(x,h) \alpha^{k-1} \int_{ \frac{\alpha x}{h}}^{\infty} u^{-k} e^{-u} \mathrm{d}u \nonumber \\
        & = c(x,h) \alpha^{k-1} \Gamma\left(-k+1, \frac{\alpha x}{h}\right).        
    \end{align}
\end{IEEEproof}
\corollary\label{cor:barthetaxh} The CDF, PDF and mean of ${\Theta}^T$ conditioned on $(X^+,H^+)=(x,h)$ are
\begin{align}
    \mathbb{P}[{\Theta}^T_{x,h}\leq \phi] &= \exp\left(-\frac{\lambda}{\mu}\left[\frac{1}{\tan \phi}-\frac{x}{h}\right]e^{-\mu h}\right), \nonumber\\
    f_{{\Theta}^T_{x,h}}(\phi) &=\frac{\lambda}{\mu \sin^2{\phi}}\exp\left(-\mu h +\frac{\lambda e^{-\mu h }}{\mu h}\left(x-\frac{h}{\tan \phi}\right)\right),\nonumber \\
    \mathbb{E}[{\Theta}^T_{x,h}]&=\int^{\arctan \frac{h}{x}}_{0}\phi f_{{\Theta}^T_{x,h}}(\phi) 
 \mathrm{d}\phi, \nonumber
\end{align}
for $0<\phi \leq \arctan \frac{h}{x}$, respectively.
\begin{IEEEproof}
    The proof of Corollary \ref{cor:barthetaxh} is similar to that of Corollary \ref{cor1}.
\end{IEEEproof}

\begin{figure}[t]
\begin{center}
\includegraphics[scale=0.20]{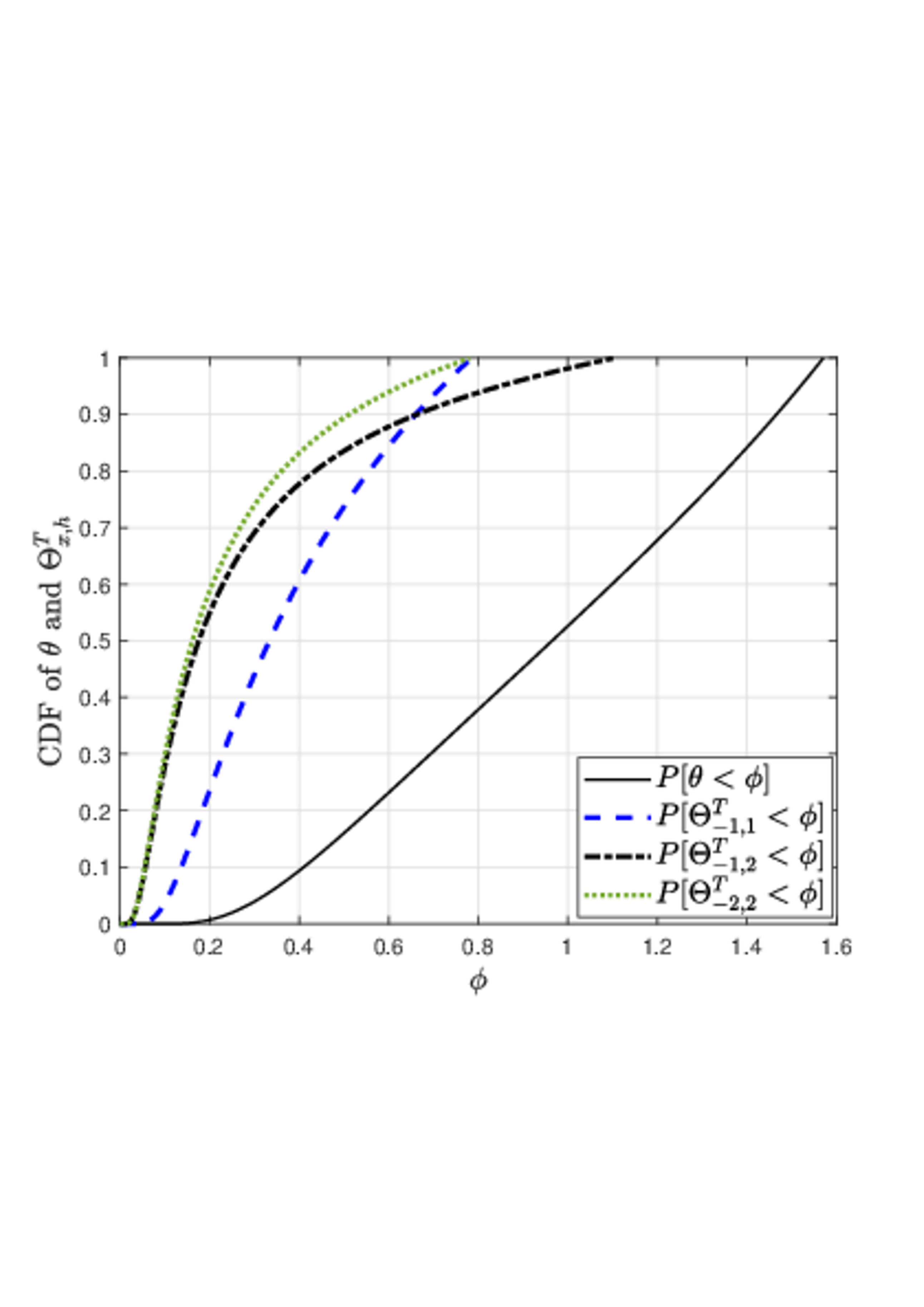}
\end{center}
\vspace{-0.5em} 
\caption{Visibility enhancement by transmissive RISs with $\lambda, \mu = 1$. } 
\label{fig:trans1} 
\end{figure}

From Corollary \ref{cor1} and Corollary \ref{cor:barthetaxh}, $\mathbb{P}[\theta\leq \phi]<\mathbb{P}[{\Theta}_{x,h}^T\leq \phi]$ for all $\phi$ in $[0, \arctan \frac{h}{x}]$. This stochastic ordering relation implies that $\mathbb{E}[\theta]>\mathbb{E}[{\Theta}^T_{x,h}]$ for all $x, h$, which shows the visibility enhancement by transmissive RISs. Figure \ref{fig:trans1} illustrates the CDFs of $\theta$ and ${\Theta}^T_{x,h}$ with $(x,h) = (1,1), (1,2)$ and $(2,2)$ for $\lambda=1$ and $\mu = 1$. In these cases, $\mathbb{E}[\theta]=0.9493$, $\mathbb{E}[{\Theta}^T_{1,1}]=0.3669$, $\mathbb{E}[{\Theta}^T_{1,2}]=0.2714$, and $\mathbb{E}[{\Theta}^T_{2,2}]=0.2249$. Under the transmissive mode, the distribution of ${\Theta_{x,h}^T}$ depends on $x$ unlike that of $\theta_{x,h}$.

\subsection{Reflective Mode}
Figure \ref{fig:sys_model_left} illustrates the visibility enhancement using the reflective property of RIS. Here, we generate an additional one-dimensional marked homogeneous PPP with intensity $\lambda$, which represents the location of buildings along the negative direction, and assume that RISs are installed on the roof of each building. The user located at the origin transmits its signal to the roof of the building, which creates the maximum blockage along the negative direction, and the RIS installed on that building reflects the signal into the sky along the positive direction.

As in the previous section, we define ${\Theta}^R_{x,h}$ and ${\Psi}^R_{x,h}$ to be the blockage angle and the visibility angle given that $(X^-,H^-)=(x,h)$, respectively, where $x<0$, as illustrated in Figure \ref{fig:sys_model_left}. In the following theorem, we first analyze the distribution of ${\Theta}^R_{x,h}$.

\theorem \label{theo4} The CDF of $\tan{\Theta}^R$ given that $(X^-,H^-)=(x,h)$ is
\begin{align}\label{eq:theo4}
&\mathbb{P}[\tan{\Theta}^R \leq t|X^-=x, H^-=h]\nonumber\\
&= \mathbb{P}[\tan {\Theta}^R_{x,h}\leq t]=\exp\left(-\frac{\rho}{ t} e^{-\mu h } e^{\mu t x}\right),
\end{align}
for $t\geq 0$ and the CDF, PDF and mean of ${\Theta}^R_{x,h}$ are
\begin{align}\label{eq:CDF_theta-}
    \mathbb{P}[{\Theta}^R_{x,h}\leq \phi] &=\exp\left(-\frac{\rho}{ \tan\phi} e^{-\mu h } e^{\mu \tan\phi x}\right),\nonumber\\
    f_{{\Theta}^R_{x,h}}(\phi)&= \exp\left(-\frac{\rho}{ \tan\phi} e^{-\mu h } e^{\mu \tan\phi x} - \mu h + \mu \tan \phi x\right)\nonumber\\
    &\times\left(\frac{\rho}{ \sin^2\phi} 
    -   \frac{\lambda x}{\sin\phi \cos \phi}
    \right),\nonumber\\
    \mathbb{E}[{\Theta}^R_{x,h}]&=\int_{0}^{\frac{\pi}{2}}\phi f_{{\Theta}^R_{x,h}}(\phi) d\phi,
\end{align}
for $0\leq \phi \leq \frac{\pi}{2}$, respectively.
\begin{IEEEproof}
See Appendix \ref{appentheo4}.
\end{IEEEproof}

\begin{figure}[t]
\begin{center}
\includegraphics[scale=0.20]{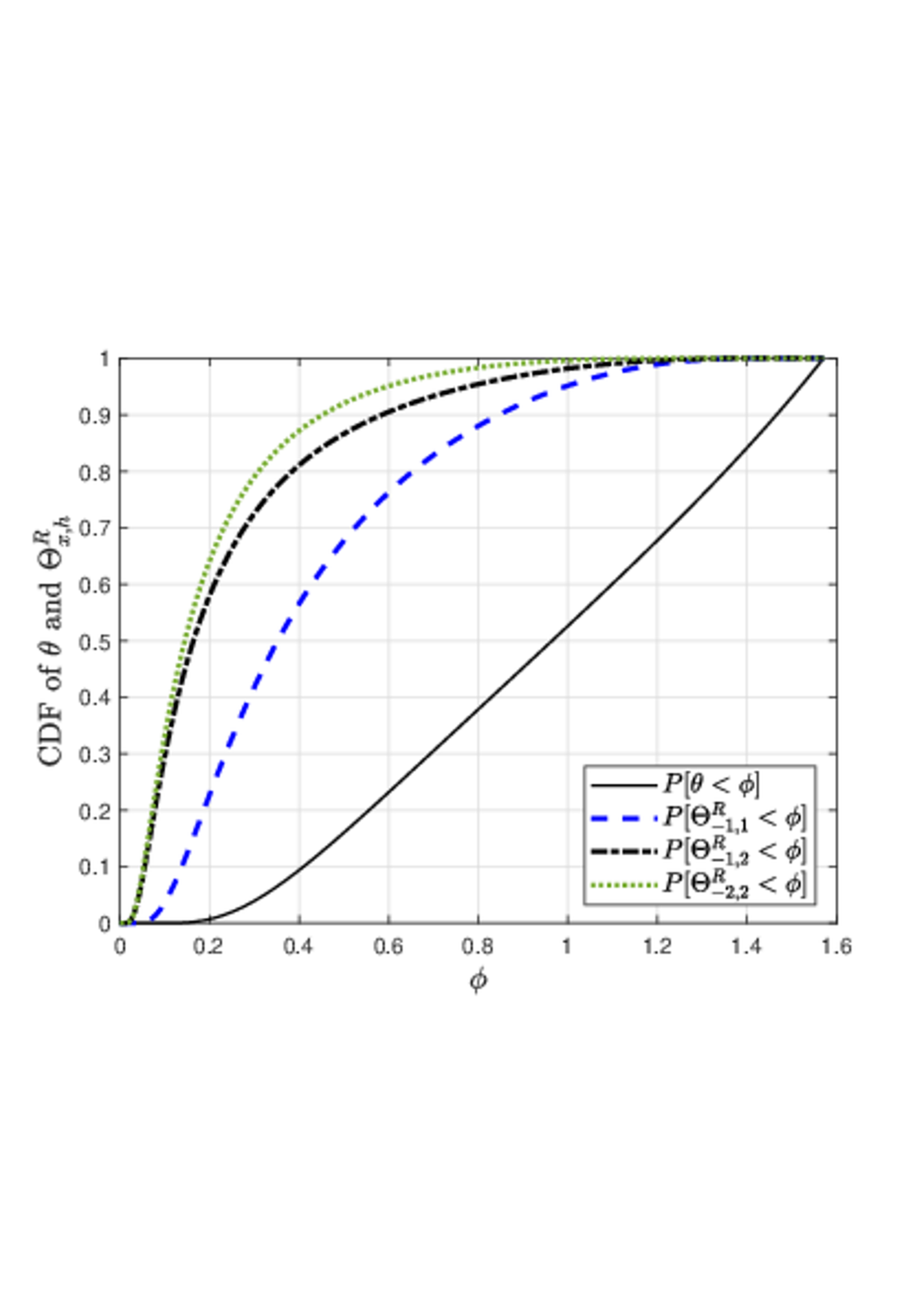}
\end{center}
\vspace{-0.5em} 
\caption{Visibility enhancement by reflective RISs with $\lambda,\mu = 1$. } 
\label{fig:refle1} 
\end{figure}

As in the transmissive mode, when $x = 0$, Theorem \ref{theo4} boils down to Theorem \ref{theo3} by plugging $x=0$, since the visibility angles observed by reflective RIS at $(0,h)$ and the user at $(0,h)$ are the same. So, $\Theta^R_{0,h}\stackrel{d}{=}\theta_{0,h}$.

Since $x<0$ in the reflective mode, $\mathbb{P}[\theta\leq \phi] < \mathbb{P}[{\Theta}^R_{x,h}\leq\phi]$, and this relation leads to $\mathbb{E}[\theta]>\mathbb{E}[{\Theta}^R_{x,h}]$ for all $x,h$. This shows the visibility enhancement by RISs under the reflective mode. Figure \ref{fig:refle1} shows the stochastic ordering relation of CDFs of $\theta$ and ${\Theta}^R_{x,h}$ for $(x,h)=(-1,1), (-1,2)$ and $(-2,2)$ where $\lambda = 1$, $\mu=1$. In these cases, $\mathbb{E}[\theta] = 0.9493$, $\mathbb{E}[{\Theta}^R_{-1,1}] = 0.4272$, $\mathbb{E}[{\Theta}^R_{-1,2}] = 0.2505$, and $\mathbb{E}[{\Theta}^R_{-2,2}] = 0.2068$. As in the transmissive mode, the distribution of ${\Theta}^R_{x,h}$ depends on $x$.

\section{Metrics, Numerical Experiments and Discussion}\label{sec: Section5}

This section focuses on the visibility enhancement obtained by RISs based on numerical experiments. We consider the visibility gain in two domains, which are 1) the angular domain and 2) the linear domain. We first study the relation between the network environment and the enhancement of visibility angle, $\psi$, under both the transmissive and reflective RIS modes. We also compare these two RIS modes regarding visibility enhancement in the angular domain. Then, we will consider how many aerial nodes can be observed by the transmissive RISs in a two-hop manner and which are not observed directly by the ground user.  

\subsection{Visibility Gain in the Angular Domain}

\begin{figure}[t]
\begin{center}
\includegraphics[scale=0.20]{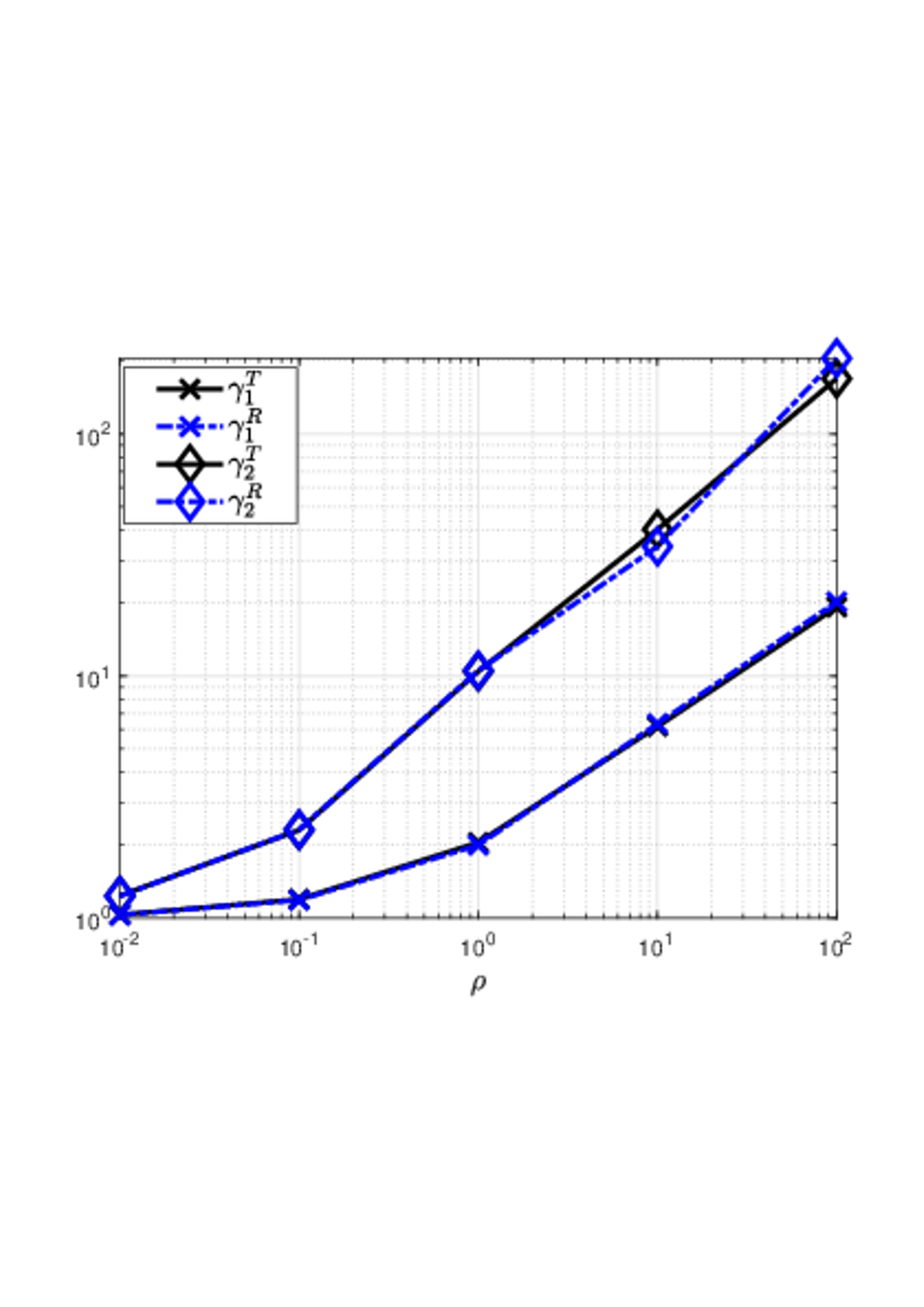}
\end{center}
\vspace{-0.5em} 
\caption{Visibility enhancement with respect to $\rho$. } 
\label{fig:metric1} 
\end{figure}

In order to quantify the visibility enhancement by RISs, we define 
\begin{align}
&\gamma_1^T \triangleq \frac{\mathbb{E}[\Psi^T]}{\mathbb{E}[\psi]}, ~~\gamma_1^R \triangleq \frac{\mathbb{E}[\Psi^R]}{\mathbb{E}[\psi]},\nonumber\\
&\gamma_2^T \triangleq \mathbb{E}\left[ \frac{\Psi^T}{\psi}\right], \gamma_2^R \triangleq \mathbb{E}\left[ \frac{\Psi^R}{\psi}\right],
\end{align}
which measure how much the visibility angle is increased by using the RISs. Figure \ref{fig:metric1} depicts the visibility enhancement with $\mu=1$ by varying $\lambda$. We can observe that 1) the visibility gain by RISs is more effective when the outdoor environment is denser, and 2) the gains between the transmissive mode and the reflective mode are not very different. Intuitively, when the density and height of buildings and the height are very low, the visibility angle observed by the user, $\psi$, is almost $\frac{\pi}{2}$. 

Further, we can observe that $\gamma_1^T$ and $\gamma_1^R$ are almost identical and $\gamma_2^T$ and $\gamma_2^R$ are almost the same for all $\rho$, which implies that we can expect the same visibility gains from transmissive RISs and reflective RISs for any outdoor environment of the M/M type.

\subsection{Linear Visibility Gain}\label{subsec:lvg}

We investigate how many additional non-terrestrial nodes can be connected through RISs. {\color{black}
We focus on the transmissive mode since the angular gains obtained by RIS between the two modes are not significant.} We assume that the ground user is located at $(0,0)$ and transmissive RISs are installed on the roofs of buildings. We assume that aerial nodes are deployed at the same altitude. 

\begin{figure}[t]
\begin{center}
\includegraphics[scale=0.27]{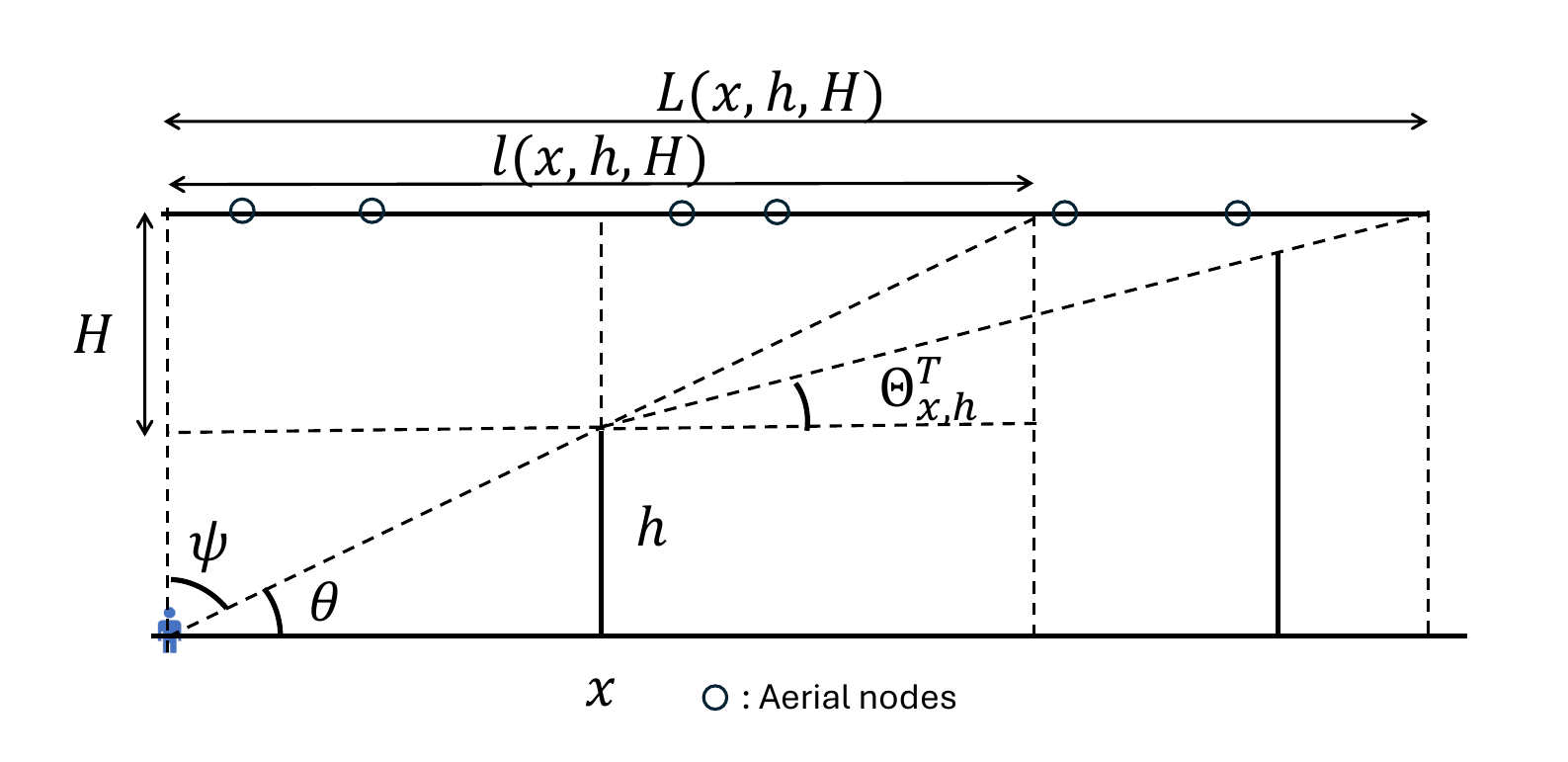}
\end{center}
\vspace{-0.5em} 
\caption{Extension of visible regions by transmissive RISs.} 
\label{fig:aerial_nodes} 
\end{figure}

{ \footnotesize
\begin{table}
\begin{center}
\begin{tabular}{|c || c|} 
  \hline
    \textbf{Metrics} & \textbf{Definition}\\
 \hline
    \multirow{2}{*}{$l(x,h,H)$} & Length of visible region at altitude $h+H$\\
    &without RIS conditioned on $(X^+,H^+)=(x,h)$\\
 \hline
   \multirow{2}{*}{$L(x,h,H)$} & Length of visible region at altitude $h+H$\\
    &with RIS conditioned on $(X^+,H^+)=(x,h)$\\
 \hline
    \multirow{2}{*}{$l(H)$}   & Length of visible region without RIS obtained by \\
    & deconditioning $l(x,h,H)$ with respect to $j(x,h)$\\
 \hline 
 \multirow{2}{*}{$L(H)$}   & Length of visible region with RIS obtained by \\
     & deconditioning $l(x,h,H)$ with respect to $j(x,h)$\\
 \hline
 \multirow{2}{*}{$\tau_H(x,h)$}& Conditional connectivity probability\\
 &when $(X^+,H^+)=(x,h)$\\
 \hline
 $\tau_H$& Conditional connectivity probability\\
\hline
\end{tabular}
\end{center}
\vspace{-0.5em}
\caption{Definitions of metrics to measure linear visibility gain in Sections \ref{subsec:lvg} and \ref{subsec:cs}.}
\end{table}
}

\subsubsection{Extension of the Length of the Visible Region}

Let $x$ and $h$ be the location and the height of the blocking building along the positive direction. As depicted in Figure \ref{fig:aerial_nodes}, we define $L(x,h,H)$ and $l(x,h,H)$ to be the visible regions at altitude $h+H$, with and without RISs, conditioned on $(X^+,H^+)=(x,h)$, where $H$ is a fixed constant. Then, 
\begin{align}
    |l(x,h,H)| & = x + H \frac{x}{h},\\
    |L(x,h,H)| & = x+\frac{H}{\tan \Theta^T_{x,h}},
\end{align}
and $|L(x,h,H)|-|l(x,h,H)|$ is the extended length of the visibility region by transmissive RISs.

The expectation of the length of $L(x,h,H)$ is
\begin{align}
    \mathbb{E}[|L(x,h,H)|]&=x+\int_0^{\frac{h}{x}} \frac{H}{t}  f_{\tan {\Theta}_{x,h}^T}(t) \mathrm{d}t\nonumber\\
    &=
    x+\frac{e^{h\mu}h + x \rho}{h\rho}H,
\end{align}
where $f_{\tan {\Theta}_{x,h}^T}(t)$ is in \eqref{eq:trans_tan}. So, conditioned on $(X^+,H^+)=(x,h)$, RISs extend the visible length at altitude $H+h$ as
\begin{align}
\mathbb{E}[|L(x,h,H)|-|l(x,h,H)|]=\frac{e^{h\mu}}{\rho}H.
\end{align}

Let $l(h)$ and $L(h)$ be the lengths of the visible region obtained by deconditioning $l(x,h,H)$ and $L(x,h,H)$ with respect to $j(x,h)$ given in Theorem \ref{theo2}. We get that the expectations of $l$ and $L$ are
\begin{align}
    \mathbb{E}[l] &= \int_{0}^{\infty}\int_{0}^{\infty} \left(x+ H\frac{x}{h}\right)j(x,h) \mathrm{d}x \mathrm{d}h \nonumber\\
    &= \frac{2+H\mu}{\lambda},\nonumber\\
    \mathbb{E}[L] & = \mathbb{E}[l]+\int_{0}^{\infty}\int_{0}^{\infty} \frac{e^{h\mu} }{\rho}H  j(x,h)   \mathrm{d}x \mathrm{d}h \nonumber\\
    &=  \frac{2+H\mu}{\lambda} 
+\int_0^{\infty} \frac{h\mu^2}{\rho} \mathrm{d}h = \infty .\nonumber
\end{align}

{\color{black}Thus, if aerial nodes are distributed as a stationary point process with intensity $\nu$, the user is connected to $\frac{2+H\mu}{\lambda}\nu$ aerial nodes without RIS in an average sense, whereas the user can connect to an infinite mean number of aerial nodes by using RISs in a two-hop manner.}

\subsubsection{Conditional Coverage Probability by RISs}

Another important metric is the conditional coverage probability, namely the probability that at least one aerial node is connected to the ground user with transmissive RISs when no aerial nodes are connected to the ground user without RIS. Let $\Phi_u$ be the location of HAPs, which are assumed to be distributed as a homogeneous PPP with intensity $\nu$ at altitude $h+H$. Let $\tau_H(x,h)$ be the conditional connectivity probability given that $(X^+,H^+)=(x,h)$. Then, 
\begin{align}
	&\tau_H(x,h)= \mathbb{P}[\Phi_u({L}(x,h,H))>0|\Phi_u(\ell(x,h,H))=0]\nonumber\\
	& = \frac{\mathbb{P}[\Phi_u({L}(x,h,H))>0,\Phi_u(\ell(x,h,H))=0]}{\mathbb{P}[\Phi_u(\ell(x,h,H))=0]}\nonumber\\
	& = \frac{\mathbb{P}[\Phi_u({L}(x,h,H)\setminus \ell(x,h,H))>0] \mathbb{P}[\Phi_u(\ell(x,h,H))=0]}{\mathbb{P}[\Phi_u(\ell(x,h,H))=0]}\nonumber\\
    & \stackrel{(a)}{=} 1-\mathbb{E}\left[\exp\left(-\nu \left(x+\frac{H}{\tan \Theta^T_{x,h}}-x-H\frac{x}{h}\right)\right)\right]\nonumber\\
    & = 1- \int_{0}^{\frac{h}{x}}e^{-\nu H \left(\frac{1}{t}-\frac{x}{h}\right)} f_{\tan {\Theta}_{x,h}^T}(t) \mathrm{d}t\nonumber\\
    & = 1-\frac{\rho}{e^{h\mu}H\nu + \rho},\nonumber
\end{align}
where $(a)$ comes from the probability generating functional (PGFL) of the PPP \cite{baccelli2009stochastic}. 

\begin{figure}[t]
\begin{center}
\includegraphics[scale=0.20]{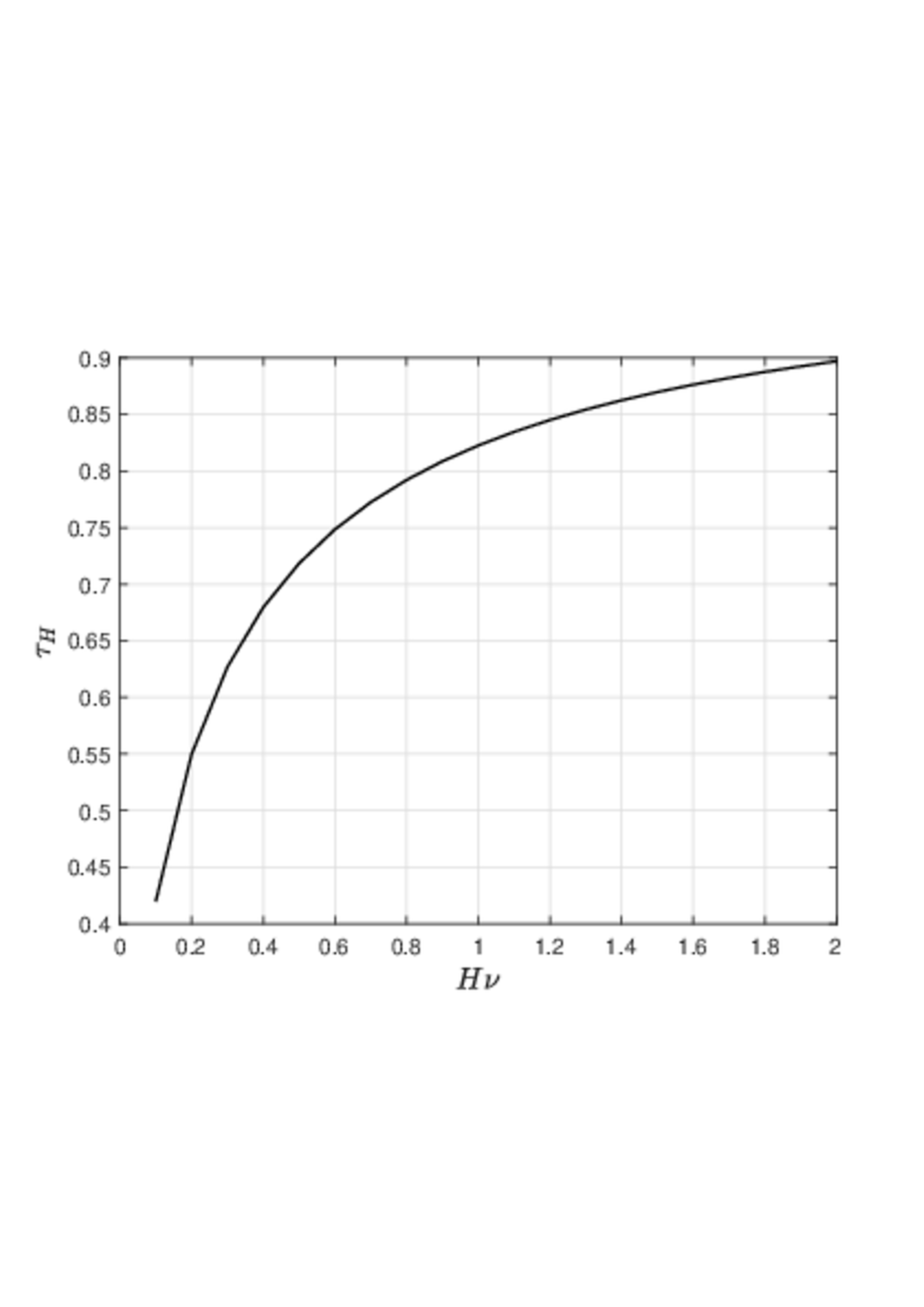}
\end{center}
\vspace{-0.5em} 
\caption{$\tau_H$ with respect to $H\nu$ for $\lambda,\mu = 1$.} 
\label{fig:visible_nodes} 
\end{figure}

Further, we define $\tau_H$ as the unconditioned conditional connectivity probability with respect to $j(x,h)$ in Theorem \ref{theo2}. Then, we obtain that
\begin{align}\label{eq:tau}
    &\tau_H = \frac{H\nu}{6\rho}\Bigg(\pi^2+6\log\left(\frac{\rho}{H\nu}\right)\log\left(1+\frac{\rho}{H\nu}\right)\nonumber\\
    &-3\left(\log\left(1+\frac{\rho}{H\nu}\right)\right)^2-6Li_2\left(\frac{H\nu}{H\nu +\rho}\right)\Bigg),
\end{align}
where $Li_n(z)$ is the polylogarithm function  
\begin{align}
    Li_n(z)=\sum_{k=1}^{\infty}\frac{z^k}{k^n}.
\end{align}
So, $\tau_H$ is a function of $\nu H$ and illustrated in Figure \ref{fig:visible_nodes}.

\subsection{Case study}\label{subsec:cs}

{ \footnotesize
\begin{table}
\begin{center}
\begin{tabular}{|c ||c| c|c|} 
 \hline
   \multirow{2}{*}{} & \textbf{Case 1} & \textbf{Case 2} & \textbf{Case 3}\\
   & \textbf{(Dense urban)} & \textbf{(Urban)} & \textbf{(Suburban)}\\
 \hline
    $\lambda (m^{-1})$& $0.012$ & $0.007$ & $0.001$\\
 \hline
    $\mu (m^{-1})$ & $0.02$ & $0.02$ & $0.02$\\
 \hline
  $\mathbb{E}[\theta]$ (rad)  & $0.7732$  & $0.5935$  & $0.1695$ \\
 \hline
 $\mathbb{E}[\theta^T]$ (rad)& $0.1956$ & $0.1256$  & $0.0195$ \\
 \hline
 $\mathbb{E}[\theta^R]$ (rad)& $0.2214$ & $0.1453$  & $0.0201$ \\
\hline
 {\color{black}$\mathbb{E}[l]$ (HAP) (m)} & $1.68\times 10^4$& $2.89 \times 10^4$& $2.02 \times 10^5$\\
\hline
 $\mathbb{E}[l]$ (Satellite) (m) & $8.34\times 10^5$& $1.43 \times 10^6$& $1.00 \times 10^7$\\
\hline
 {\color{black}$\tau_H$ (HAP)} & $0.797838$&$0.864512$&$0.976052$ \\
\hline
 $\tau_H$ (Satellite)& $0.893742$&$0.935147$ &$0.989409$\\
\hline
\end{tabular}
\end{center}
\vspace{-0.5em}
\caption{{\color{black}Numerical values for visibility in the three cases.}}
\end{table}
}

{\color{black}
We provide three examples to show the enhancement of visibility. These examples represent 1) dense urban, 2) urban and 3) suburban or rural areas. We set $\rho = 0.6 (\lambda = 0.012 (m^{-1},\mu =0.02 (m^{-1}), 0.35 (\lambda = 0.007 (m^{-1},\mu =0.02 (m^{-1}), 0.1 (\lambda = 0.001 (m^{-1},\mu =0.02 (m^{-1})$ for these examples as explained in Figure \ref{fig:response1_3_2}. For these cases, we consider two scenarios: 1) HAP with $\nu = 5\times 10^{-5} m^{-1}$ and $H = 10km$, and 2) Satellite with $\nu =2.3163\times 10^{-6} m^{-1}$ and $H = 500 km$, and Table III shows the metrics of sky visibility for these three examples. 

{\color{black}
From this table, we can observe that the mean values of $\theta$, $\theta^T$ and $\theta^R$ increase as $\rho(=\frac{\lambda}{\mu})$ increases. Further, the mean length of the visibility region without RISs, $\mathbb{E}[l]$, increases as $\rho$ decreases. We also provide numerical results on $\tau_H$, the probability that there exists at least one observable aerial node thanks to RISs when the user cannot connect to any aerial nodes without RIS, for these three cases. For the dense urban case (Case 1), $\tau_H$ for the HAP case is $0.797838$. This can be interpreted as the fact that the probability that the ground user can observe at least one HAP via a RIS in a two-hop LOS manner is around $80\%$ when he/she cannot observe any HAP without RIS. $\tau_H$ also increases as $\rho$ decreases. This comes from the fact that the sky visibility at the RIS on the blocking building is still better when $\rho$ is small under our M/M model.
}

}

\section{Conclusion}\label{sec: Section6}

In this paper, we proposed a novel framework to analyze the user's visibility into the sky in an outdoor environment using the theory of point processes and stochastic geometry. Based on that model, we presented the distributions of the blockage and visibility angles as a function of network parameters. We investigated the visibility enhancement when transmissive RISs and reflective RISs installed on the roofs of buildings were used as relay nodes in a two-hop manner. Finally, through the analytical results, we analyzed how much visibility improved thanks to RISs in terms of angular and linear visibility gains.

\appendices
\section{Proof of Theorem \ref{theo1}}\label{appen1}
The CDF of $\tan\theta$ is
\begin{eqnarray}
 &\mathbb{P}\left[\tan \theta \leq t\right] &= \mathbb{P}\left[\max f(x_i,h_i) \leq t\right]\nonumber\\
& &  =\mathbb{E}\left[\mathbbm{1}(\max f(x_i,h_i)\leq t)\right]\nonumber\\
&  & =\mathbb{E}\left[\prod_{x_i\in\Phi^+} \mathbbm{1}(f(x_i,h_i)\leq t)\right]\nonumber\\
    && = \mathbb{E}\left[\exp\left(\sum_{x_i\in\Phi^+}\log (\mathbbm{1}(f(x_i,h_i)\leq t))\right)\right]\nonumber\\
    && \stackrel{(a)}{=} \exp\left(-\lambda\int_{0}^{\infty}1-\mathbb{E}_h[\mathbbm{1}(f(x,h)\leq t)]\mathrm{d}x\right)\nonumber\\
    &&=  \exp\left(-\lambda\int_{0}^{\infty}1-\mathbb{P}[h \leq xt]\mathrm{d}x\right)\nonumber\\
    && =  \exp\left(-\lambda\int_{0}^{\infty}e^{-(\mu x t)^k}\mathrm{d}x\right)\nonumber\\
    && =  \exp\left(-\frac{\lambda}{\mu t}\Gamma\left(1+\frac{1}{k}\right)\right), \nonumber\\
    && = \exp\left(-\frac{\rho}{t}\Gamma\left(1+\frac{1}{k}\right)\right),\nonumber
\end{eqnarray}
where $\mathbbm{1}(\cdot)$ is the indicator function, and $(a)$ comes from the probability generating functional (PGFL) of PPP. \cite{baccelli2009stochastic}

\section{Proof of Variant \ref{ex2}}\label{appenvar2}
    Conditioned on $x_1=u\in[0,\frac{1}{\lambda})$, the CDF of $\tan {\theta}$ is  
    \begin{eqnarray}
    &&\mathbb{P}[\tan {\theta} \leq t|x_1=u]\nonumber\\
    &&= \mathbb{P}[\tan \theta_1 \leq t, \tan \theta_2 \leq t, ...|x_1=u]\nonumber\\
    &&=\mathbb{P} \left[h_1\leq tu,h_2 \leq t\left(u+\frac{1}{\lambda}\right), ...\right]\nonumber\\
    &&= \prod_{i=1}^{\infty}\mathbb{P}\left[h_i\leq t\left(u+\frac{i-1}{\lambda}\right)\right]\nonumber\\
    &&= \prod_{i=1}^{\infty}\left(1-\exp\left(-\mu\left(u+\frac{i-1}{\lambda}\right)t\right)\right),\nonumber
    \end{eqnarray}
    since the random variables $\{h_i\}$ are independent exponential random variables. By integrating over $u$, we can obtain the CDF of $\tan{\theta}$ as \eqref{eq6}. The CDF of ${\theta}$ can be obtained as in the proof of Theorem \ref{theo1}.

\section{Proof of Theorem \ref{theo3}}\label{appentheo3}

When the user is located at $(x,h)$, a building located at $x_i$ with a height $h_i$ for $x < x_i$ blocks the user's visibility depending on the value of $\frac{h_i-h}{x_i-x}$. For all $t$, we have 
    \begin{eqnarray}
    &&\mathbb{P}[\tan {{\theta}}_{x,h}\leq t]= \mathbb{P}\left[\max_{x_i>x} \frac{h_i-h}{x_i-x}\leq t \right]\nonumber\\
    &&= \mathbb{E}\left[\prod_{x_i\in\Phi\cap (x,\infty)}\mathbbm{1}\left(\frac{h_i-h}{x_i-x}\leq t\right)\right]\nonumber\\    
    &&= \mathbb{E}\left[\exp\left(\sum_{x_i\in\Phi\cap (x,\infty)}\log\mathbbm{1}\left(\frac{h_i-h}{x_i-x}\leq t\right)\right)\right]\nonumber\\    
    &&= \exp\left(-\lambda\int_{x}^{\infty} 1-\mathbb{E}_v \left[\mathbbm{1}\left(\frac{v-h}{u-x}\leq t\right)\right]\mathrm{d}u\right)\nonumber\\
    &&= \exp\left(-\lambda\int_{x}^{\infty} 1-\mathbb{P} \left[\frac{v-h}{u-x}\leq t\right]\mathrm{d}u\right)\nonumber\\
    &&= \exp\left(-\lambda\int_{x}^{\infty} 1-\mathbb{P} \left[v \leq t(u-x)+h\right]\mathrm{d}u\right)\nonumber\\
    &&= \exp\left(-\lambda e^{-\mu h}\int_{0}^{\infty} e^{-\mu tz} \mathrm{d}z\right)\nonumber\\
    &&= \exp\left(-\frac{\lambda}{\mu t} e^{-\mu h}\right).\nonumber
\end{eqnarray}
The CDF and the PDF of ${\theta}_{x,h}$ can be obtained as the proof of Theorem \ref{theo1}.

\section{Proof of Theorem \ref{theo2}}\label{appentheo2}
Let $\gamma_{i,\lambda}(x)=\frac{\lambda^i x^{i-1}e^{-\lambda x}}{(i-1)!}$ denote
the density of Erlang distribution with shape parameter $i$ and rate $\lambda$.
The joint distribution of $X^+$ and $H^+$ is 
\begin{align}
	&\mathbb{P}[X^+\in [x,x+\mathrm{d}x],H^+\in[h,h+\mathrm{d}h] ]\nonumber\\
	&=\sum_{i\geq 1}\mathbb{P}\Bigg[x_i\in [x,x+\mathrm{d}x], \nonumber\\
 &~~~~~~~~~~h_i\in[h,h+\mathrm{d}h],\bigcap_{j\neq i} \frac{h_j}{x_j}<\frac{h_i}{x_i}\Bigg]\nonumber\\
	&=\sum_{i\geq 1}\mathbb{P}\Bigg[\bigcap_{j\neq i} \frac{h_j}{x_j}<\frac{h_i}{x_i} \bigg| x_i\in [x,x+\mathrm{d}x],\nonumber\\ & \hspace{1cm}h_i\in[h,h+\mathrm{d}h] \Bigg]\mathbb{P}[x_i\in [x,x+\mathrm{d}x]]\nonumber\\
    &\times\mathbb{P}[h_i\in[h,h+\mathrm{d}h]]\mathrm{d}x\mathrm{d}h.\nonumber
\end{align}
That is
\begin{align}
	&\mathbb{P}[X^+\in [x,x+\mathrm{d}x],H^+\in[h,h+\mathrm{d}h] ]\nonumber\\
    &\stackrel{(a)}{=} \sum_{i\geq 1}\left[\frac{1}{x}\int_{0}^{x}\mathbb{P}_g\left[\frac{g}{y}<\frac{h}{x}\right] dy\right]^{i-1}\nonumber\\
    &\times \exp\left(-\lambda\int_{x}^{\infty}e^{-\mu v\frac{h}{x}}\mathrm{d}v\right) \gamma_{i,\lambda}(x) \mu e^{-\mu h}\nonumber\\
	&= \sum_{i\geq 1}\gamma_{i,\lambda}(x)\mu e^{-\mu h}\left[\frac{1}{x}\int_{0}^{x}\left(1-e^{-\frac{\mu hy}{x}}\right) dy\right]^{i-1}\nonumber\\
    &\times \exp\left(-\frac{\lambda x}{\mu h}e^{-\mu h}\right)\mathrm{d}x\mathrm{d}h\nonumber\\
    &=\sum_{i\geq 1}\gamma_{i,\lambda}(x)\mu e^{-\mu h} \left(1-\frac{1-e^{-\mu h}}{\mu h}\right)^{i-1} \nonumber\\
    &\times \exp\left(-\frac{\lambda x}{\mu h}e^{-\mu h}\right)\mathrm{d}x\mathrm{d}h\nonumber\\
	&=\mu e^{-\mu h}\exp\left(-\frac{\lambda x}{\mu h}e^{-\mu h}\right)\nonumber\\
 &\times\sum_{i\geq 1}\frac{\lambda^i x^{i-1}e^{-\lambda x}}{(i-1)!} \left(1-\frac{1-e^{-\mu h}}{\mu h}\right)^{i-1}\mathrm{d}x\mathrm{d}h \nonumber\\
	&=\lambda\mu e^{-\mu h}\exp\left(-\frac{\lambda x}{\mu h}e^{-\mu h}-\lambda x\right)\nonumber\\
 &\times\sum_{i\geq 1}\frac{ 1}{(i-1)!} \left(\lambda x\left(1-\frac{1-e^{-\mu h}}{\mu h}\right)\right)^{i-1} \mathrm{d}x\mathrm{d}h\nonumber\\
	&\stackrel{(b)}{=}\lambda\mu e^{-\mu h}\exp\left(-\frac{\lambda x}{\mu h}e^{-\mu h}-\lambda x\right) e^{\lambda x \left(1-\frac{1-e^{-\mu h}}{\mu h}\right)}\mathrm{d}x\mathrm{d}h \nonumber\\
	&=\lambda\mu \exp\left(-\mu h-\frac{\lambda x}{\mu h}\right)\mathrm{d}x\mathrm{d}h. \nonumber
\end{align}
In $(a)$, we used the fact that the density of $x_i$ is $\gamma_{i,\lambda}$ and that of $h_i$ is $\mu e^{-\mu h}$
to get the last two terms. 
The first term is the conditional probability of the event $\bigcap_{1\le j < i} \frac{h_j}{x_j}<\frac{h}{x}$.
To get this expression, we use the well-known property of the PPP that conditionally on the fact that $x_i=x$,
the $i-1$ (unordered) points of the PPP in the $[0,x]$ interval are independently and uniformly distributed.
In this first term, $g$ is an independent exponential random variable with parameter $\mu$.
The second term is the conditional probability that no buildings after $x_i$ affect
the blockage observed by the user at $0$. The relation in
$(b)$ comes from the Taylor series expansion of the exponential function.


\section{Proof of Theorem \ref{theo:vis_en_tr_tan}}\label{appentheo_vis_en_tr_tan}

Let ${\Theta}^T$ denote the blockage angle seen from $(X^+,H^+)$ in the transmissive mode.
By definition $\mathbb{P}[\tan {\Theta}^T_{x,h} \leq t]  =\mathbb{P}[\tan {\Theta}^T \leq t | X^+ = x, H^+ = h]$
and by arguments similar to those used in the last proof,
\begin{align}\label{eq:theo_vis_en_tr_tan2} 
    &\mathbb{P}[\tan {\Theta}^T \leq t \big| X^+ =x , H^+ = h]\nonumber\\&= \sum_{i=1}^{\infty}\mathbb{P}\Bigg[ \bigcap_{j< i} \frac{h_j}{x_j}<\frac{h_i}{x_i},\bigcap_{j>i} \frac{h_j-h_i}{x_j-x_i}<t\bigg|x_i=x,h_i=h\Bigg] \nonumber\\
    &=\sum_{i \geq 1}\left[ \frac{1}{x}\int_{0}^{x}\mathbb{P}_g\left[\frac{g}{y}<\frac{h}{x} \right]dy\right]^{i-1}\nonumber\\
    &\times \exp\left(-\lambda\int_{x}^{\infty}e^{-\mu (t(v-x)+h)}\mathrm{d}v\right).
\end{align}
Hence,
\begin{align}\label{eq:theo_vis_en_tr_tan22} 
	&\mathbb{P}[\tan {\Theta}^T \leq t, X^+ \in [x,x+\mathrm{d}x] , H^+ \in [h,h+\mathrm{d}h]]
	\nonumber\\&= \sum_{i=1}^{\infty} \gamma_i(x) e^{-\mu h} \left(1-\frac{1-e^{-\mu h}}{\mu h}\right)^{i-1}\nonumber\\&\times \exp\left(-\lambda\int_{x}^{\infty}e^{-\mu ( t (v-x) + h)} \mathrm{d}v\right) \mathrm{d}x\mathrm{d}h \nonumber\\
    &= \sum_{i=1}^{\infty} \gamma_i(x) e^{-\mu h} \left(1-\frac{1-e^{-\mu h}}{\mu h}\right)^{i-1}e^{-\frac{\lambda}{\mu t} e^{-\mu h}}\mathrm{d}x\mathrm{d}h\nonumber\\
    & = \lambda \mu \exp\left(-\mu h -\frac{\lambda}{\mu t}e^{-\mu h} -\frac{\lambda x}{\mu h}(1-e^{-\mu h})\right)\mathrm{d}x\mathrm{d}h.
\end{align}
Equation \eqref{eq:vis_en_tr_tan} follows from this and
from the fact that (Theorem \ref{theo2})
\begin{align}\label{eq:theo_vis_en_tr_tan3} 
	&\mathbb{P}[ X^+ \in[x,x+\mathrm{d}x], H^+ \in [h+\mathrm{d}h]]\nonumber\\&= j(x,h) \mathrm{d}x\mathrm{d}h\nonumber\\&
	= \lambda\mu e^{-\mu h -\frac{\lambda x}{\mu h}} \mathrm{d}x\mathrm{d}h.
\end{align}

\section{Proof of Theorem \ref{theo4}}\label{appentheo4}

The CDF of $\tan{\Theta}^R_{x,h}$ is
\begin{align}
    &\mathbb{P}[\tan {\Theta}^R_{x,h}\leq t]=\mathbb{P}\left[\max_{x_i\in\Phi} \frac{h_i-h}{x_i-x}\leq t\right]\nonumber\\
    &=\mathbb{E}\left[\prod_{x_i\in\Phi}\mathbbm{1}\left(\frac{h_i-h}{x_i-x}\leq t \right)\right]\nonumber\\
    &=\mathbb{E}\left[\exp\left(\sum_{x_i\in\Phi}\log\left(\mathbbm{1}\left(\frac{h_i-h}{x_i-x}\leq t \right)\right)\right)\right]\nonumber\\
	&=\exp\left(-\lambda\int_{0}^{\infty}1-\mathbb{E}_v\left[\mathbbm{1}\left(\frac{v-h}{u-x}\leq t \right)\right]\mathrm{d}u\right)\nonumber\\
	&=\exp\left(-\lambda\int_{0}^{\infty}1-\mathbb{P}\left[\frac{v-h}{u-x}\leq t \right]\mathrm{d}u\right)\nonumber\\
	&=\exp\left(-\lambda\int_{0}^{\infty}1-\mathbb{P}\left[v\leq t (u-x)+h \right]\mathrm{d}u\right) \nonumber\\
	&=\exp\left(-\lambda\int_{0}^{\infty}e^{-\mu (t (u-x)+h)}\mathrm{d}x\right)\nonumber\\
	&=\exp\left(-\lambda e^{-\mu h } \int_{-x}^{\infty}e^{-\mu t z}\mathrm{d}z\right)\nonumber\\
	&=\exp\left(-\frac{\lambda}{\mu t} e^{-\mu h } e^{\mu t x}\right).\nonumber
\end{align}
The CDF and the PDF of ${\Theta}^R_{x,h}$ can be obtained as in the proof of Theorem \ref{theo1}.
\section*{Future Work}
The next natural step along this line of research will
be the analysis of the full 3-D model alluded to in the paper.
Several other questions should also be answered within
the proposed framework. Some in the direct continuation of the
basic model like two or more RIS-relaying hops, others
aiming at extending the analysis in new and complementary
directions like the sensing/signalling question discussed in the
introduction.

\section*{Acknowledgements}
The work of J. Lee was supported by the National Research Foundation of Korea(NRF) grant funded by the Korea government(MSIT) (No. 2022R1G1A1008552).

The work of F. Baccelli was supported
by the European Research Council project titled NEMO under grant ERC 788851,
by the Horizon Europe project titled INSTINCT under grant SNS 101139161,
and by the French National Agency for Research project titled
{\em France 2030 PEPR réseaux du Futur} under grant ANR-22-PEFT-0010.

\bibliographystyle{ieeetran}
\bibliography{referenceBibs}

\end{document}